\documentclass[11pt]{amsart}
\usepackage{amsbsy,amsfonts}
\usepackage{latexsym,eucal,amsmath,amsthm}
\usepackage{amssymb}
\newcommand{\db}{\mathbb }

\newcommand{\dbR}{{\db R}}

\newcommand{\dbN}{{\db N}}
\newcommand{\dbT}{{\db T}}
\newcommand{\F}{\mbox{{\em F}}}

\newcommand{\om}{\omega}

\newcommand{\x}{\xi}

\newcommand{\lam}{\lambda}

\newcommand{\la}{\lambda}

\def\R{{\db R}}
\newcommand{\half}{{\frac{1}{2}}}

\newtheorem{theorem}{Theorem}
\theoremstyle{remark}
\newtheorem{remark}{Remark}
\theoremstyle{proposition}
\newtheorem{proposition}{Proposition}
\theoremstyle{lemma}
\newtheorem{lemma}{Lemma}
\theoremstyle{corollary}
\newtheorem{corollary}{Corollary}

\numberwithin{equation}{section}
\numberwithin{lemma}{section}
\numberwithin{remark}{section}
\numberwithin{proposition}{section}

\theoremstyle{definition}
     \newtheorem{definition}[subsection]{Definition}

\begin{document}
\date{}
\title{Local well-posedness for dispersion generalized Benjamin-Ono
equations}
\author{J. Colliander}
\thanks{J.C. is supported in part by N.S.F. Grant DMS 0100595 and N.S.E.R.C.
Grant RGPIN 250233-03.}
\address{\small University of Toronto}

\author{C. Kenig}
\thanks{C.K. was supported in part by N.S.F Grant DMS 9500725}
\address{\small University of Chicago}

\author{G. Staffilani}
\thanks{G.S. was supported in part by N.S.F. Grant DMS 0100375, the Terman Award and a grant by the Sloan Foundation.}
\address{\small M.I.T.}


\begin{abstract}
In this paper we study  local well-posedness in the energy space
for a family of dispersive
equations that can be seen as dispersive ``interpolations'' between the KdV
and the Benjamin-Ono equation.
\end{abstract}

\maketitle
\section{Introduction}
We consider the initial value problem 
\begin{equation}
\left\{ \begin{array}{l}
\partial_t u + \partial_xD_{x}^{1+a} u +
\frac{1}{2}\partial_x(u^2) = 0, \\
 u(x,0) = u_0(x) \hspace{1.5cm} x \in \dbR, \,  t \in \dbR,
\end{array}\right.
\label{ivp}\end{equation}
where $0\leq a \leq 1$. Here $D_x^{1+a}$ is the Fourier multiplier operator with symbol
$|\xi|^{1+a}$. These equation arise as mathematical models for the
weakly nonlinear propagation of long waves in shallow  channels.
We recall that when $a=0$ the equation in
\eqref{ivp} is called the Benjamin-Ono equation, and when $a=1$
it is called the KdV equation. In the endpoint cases cases ($a=0$ and $a=1$), the  equations
have an infinite number of conserved integrals and are integrable by
the inverse scattering method \cite{AF83}, \cite{CoifWick90}. When $0<a<1$ there is no
integrability, but three integrals are still conserved
\cite{SSS}:
\begin{eqnarray}
    \label{overage} I_{1}&=&\int u(x,t)\,dx,\\
    \label{l2con} I_{2}&=&\int |u(x,t)|^{2}\,dx,\\
\label{h} I_{3}&=&\frac{1}{6}\int u^{3}\,dx+
\int |D^{\frac{1+a}{2}}_{x}u|^{2}\,dx.
\end{eqnarray}
Several papers have been published on the well-posedness for the
initial value problem \eqref{ivp}, below we only recall the most recent
results:

\noindent{\underline{$a=0:$}} 
For $a=0$, \eqref{ivp} is the Benjamin-Ono initial value problem which
is known to have global weak solutions in $L^{2}$ (
\cite{S}, \cite{GV1}, \cite{GV2} and \cite{T}). Moreover, Ponce \cite{P}
proved global well-posedness  in $H^{3/2}$ by first proving local
well-posedness via the energy method enhanced with dispersive
smoothing and globalizing this result with the next conservation law in the
hierarchy of conserved quantities for Benjamin-Ono. 

\noindent{\underline{$0 < a < 1$:}} Kenig, Ponce and Vega
\cite{KPV0} have shown that \eqref{ivp} is locally well-posed for
data in $H^s$ provided $s \geq \frac{3}{4} (2-a)$ using the energy
method enhanced with the smoothing effect. When $a \geq \frac{4}{5}$
note that $\frac{3}{4}(2-a) \geq \half + \frac{a}{2}$. Therefore, in
the regime $\frac{4}{5} \leq a$, the conservation law \eqref{h}
and the local theory from \cite{KPV0} combine to 
prove global well-posedness of \eqref{ivp}. 

\noindent{\underline{$a=1:$}} For $a=1$, \eqref{ivp} is the KdV initial value
problem and Bourgain \cite{B1} used a fixed point argument to prove 
local (and hence global)  well-posedness in $L^{2}$.
Subsequently, Kenig, Ponce and Vega \cite{KPV4}  proved, again using a
fixed point argument, local well-posedness
in $H^{s}$ for $s>-3/4$ and Colliander, Keel, Staffilani, Takaoka and  Tao
\cite{CKSTT}  extended this to a global result. Christ, Colliander
and Tao \cite{CCT02} recently established a local{\footnote{Global
    well-posedness at $s= -\frac{3}{4}$ remains an open problem.}} 
well-posedness result at the $s=-3/4$ endpoint for KdV by conjugating an
extension of the $s=\frac{1}{4}$ local theory \cite{KPV2} for the modified KdV
equation using the Miura transform (see \cite{MiuraTransform},
  \cite{MiuraKdVSurvey}).

Recently
Molinet, Saut and Tzvetkov \cite{MST} have shown that for $0 \leq a < 1$
an $H^s$ assumption alone on the initial data is insufficient for a
proof of local well-posedness via Picard iteration or fixed point
arguments no matter what space of functions of spacetime containing $C([0,T]; H^s )$
is considered for the contraction. Thus, the natural goal of proving
local well-posedness of \eqref{ivp} in the space $H^{s_*},~s_* =
\frac{1}{2} + \frac{a}{2}$ appearing in $I_3$ is not attainable via
a fixed point argument. In particular, the type of approach used 
by Bourgain \cite{B1} and Kenig, Ponce and Vega \cite{KPV2}
\cite{KPV3} \cite{KPV4} for the KdV initial value problem is not
possible for \eqref{ivp} in the range $0 \leq a < 1$..

Two paths emerge for addressing the well-posedness issues for
\eqref{ivp}. We might first choose to abandon the fixed point approach
to proving local well-posedness and try to further enhance the
classical energy method by applying it somehow at lower
regularity. Since \eqref{ivp} is known \cite{GV1}, \cite{GV2} 
to have global weak solutions in
$L^2$, the main issue in this approach is uniqueness. The enhancements
of the classical energy method using the smoothing effect appearing in
\cite{P}, \cite{KPV0} follow this path. As an alternate approach,
we might choose to abandon $H^s$ and prove local well-posedness via a
fixed point argument in
some other space of initial data. If the norms of the data used in the
proof can be shown to be finite for all time we might also prove global
well-posedness for such initial data.
This paper follows the second path described above by considering
initial data in the space $F^{s_*}$ defined below which involves
$L^2$-type integrability with respect to a spatial weight, in addition
to $H^s$ regularity.

\begin{remark}
  \label{nlgenremark}
The present paper should be distinguished from \cite{KPVonBO}. Here,
we are considering dispersive generalizations of the Benjamin-Ono and
KdV equations all with the same quadratic nonlinearity. In
\cite{KPVonBO}, Kenig, Ponce and Vega considered higher power
generalizations of the nonlinearity within the Benjamin-Ono setting. 
Note also that the higher power cases in the dispersion
generalized setting ($a \geq 0$) were already considered in \cite{KPV0}. 
\end{remark}


Before we present our result we observe that from \eqref{l2con} and
\eqref{h} it follows that
if $u$ is a solution for the IVP \eqref{ivp} with $a \geq 0$, then
\begin{equation}
     \label{enorm}\|u(t)\|_{H^{s_{*}}}\leq C\|u_{0}\|_{H^{s_{*}}}
     \end{equation}
where $s_{*}=1/2+a/2$. 
Indeed, by Sobolev and interpolation, we have the bound
\begin{equation*}
  \left| \int u^3 dx \right| \thicksim {{\| u \|}_{H^{\frac{1}{6}}}^3} \thicksim
{{\| u \|}_{L^2}^{3(1 - \frac{1}{3+3a})}} 
{{\| u \|}_{H^{{\half + \frac{a}{2}}}}^{3 (\frac{1}{3+3a})}},
\end{equation*}
which, when combined with \eqref{h}, gives \eqref{enorm}. 
If one could prove well-posedness in an
interval of time $[0,T]$, where $T=T(\|u_{0}\|_{H^{s_{a}}})$, then
thanks to \eqref{enorm} an iteration in time would give global
well-posedness. Unfortunately the standard way of obtaining
a local result is  by a fixed point theorem in Sobolev spaces
of $L^{2}$ type,
and the recent result of Molinet, Saut and Tzvetkov \cite{MST} shows
that for  $0\leq a<1$ this cannot be done.

\begin{remark}
\label{KP}
There is another equation for which Molinet, Saut and Tzvetkov
\cite{MST1} exhibit
a counterexample showing the failure of an iteration method in
$H^{s}$: the KP-I equation. But, in two recent papers (\cite{CKS} and
\cite{CKS1}), we show that one can still use a fixed point method to
obtain well-posedness results, as long as one considers weighted
Sobolev spaces together with the classical $H^{s}$ space.
In this paper, we prove analogous well-posedness
results for the IVP \eqref{ivp} with $0\leq a<1$ using a fixed point 
theorem method
on weighted Sobolev spaces\footnote{The description of these spaces
is contained in Definitions \ref{def1} and \ref{def2}.}. 
\end{remark}

We shall also
see, by example, that our basic estimates, Propositions \ref{prop1} and
\ref{prop2}, fail when $a=0$, which explains why our results hold only for the
case $0 < a<1$, (see the Appendix).

Throughout the paper we use the following notation for
the Fourier transform:
$$\F(f)(\x)=\hat f(\x)=\int e^{ix\x}f(x)\,dx,$$
and similarly for the inverse Fourier transform we write
$$\F^{-1}(g)(x)=\check g(x)=\int e^{-ix\x}g(\x)\,d\x.$$

\begin{definition}
Let
\begin{equation}\label{exp}
s_{*}=1/2+a/2.
\end{equation}
We  define the space of functions $F^{s}$  as the
completion of the Schwartz functions through the norm
\begin{equation}\label{fs}
 \|f\|_{F^{s}}=\|f\|_{H^{s}}+\|xf\|_{H^{s-2s_{*}}},
    \end{equation}
\label{def1}\end{definition}
We forecast that the space $F^{s_{*}}$ will be the space of
the initial data in the IVP \eqref{ivp}.
We also need a space that will contain the evolution under the dynamics
dictated by \eqref{ivp} of the initial
data. To define this space we start by introducing
the dispersive function $\omega(\x)=\x|\x|^{1+a}$ associated
to the IVP \eqref{ivp}. We also denote with
$\chi_{0}(z)$ an even 
smooth characteristic
function of the interval $[-1,1]$ and  with
$\chi$  an even smooth characteristic function of the
set $\{z : 1/2<|z|<2\}$. We use the notation
$\chi_{j}(z)=\chi(2^{-j}z)$ with $j\in \dbN$, whenever a dyadic
decomposition is needed.
\begin{definition}
We define the space $X_{s}^{b}$ for $s, b
\in \dbR$,  as the
closure of the Schwartz functions through the norm
\begin{equation}\label{xsb}
 \|f\|_{X^{b}_{s}}=\sum_{j\geq 0}2^{jb}\left(\int_{\dbR^{2}}
 \chi_{j}(\lambda-\om(\x))(1+|\x|)^{2s}|\hat{f}|^{2}(\lambda,\x)
 \, d\x d\lambda\right)^{1/2},
    \end{equation}
and the space $Y_{s_{0}, s_{1}}^{b},$ for $ s_{0}, s_{1}, b\in \dbR$  as
\begin{equation}\label{ys0s1b}
Y_{s_{0}, s_{1}}^{b}=\{f : xf \in X_{s_{0}}^{b}
\mbox{ and }  tf \in X_{s_{1}}^{b}\}.
    \end{equation}
Finally, for $s, b\in \dbR$,  we set
\begin{equation}\label{zs}
Z_{s}^{b}=X^{b}_{s}\cap Y^{b}_{s-2s_{*},s}.
\end{equation}
If an interval of time $[0,T]$ is fixed, then we
say that $f \in Z_{s,T}^{b}$ if there is an extension $\tilde f$
of $f$ on the whole real line such that $\tilde f \in Z_{s}^{b}$ and
$\|f\|_{Z_{s,T}^{b}}=\inf \|\tilde f\|_{Z_{s}^{b}}$, where the
infimum is taken over all  possible extensions $\tilde f$. 
(These norms are adaptations
of the spaces introduced \cite{B1} by Bourgain for KdV and NLS to the generalized 
Benjamin-Ono setting, where, in light of the examples in \cite{MST}, we are forced into 
introducing some version of the spatial weight.)
\label{def2}\end{definition}
 We will  prove in  Theorem \ref{main} that $Z_{s_{*}, T}^{1/2}$,
where $s_{*}$ is as in \eqref{exp}, is a space of evolution of the
initial data in $F^{s_{*}}$ through the IVP \eqref{ivp}.

We are now ready to state the main result of this paper.
\begin{theorem}\label{main}
Assume that $a\in (0,1)$. Then for any $u_{0}\in F^{s}$ with $s\geq 
s_{*}=1/2+a/2$, there exist
 $T=T(\|u_{0}\|_{F^{s_{*}}})$ and a unique solution $u$ for the IVP
 \eqref{ivp}, such that $u\in Z_{s, T}^{1/2}\cap C([0,T],F^{s})$.
 Moreover given any $T'\in(0,T)$, there exists a neighborhood
 $U_{u_{0}}$ of $u_{0}\in F^{s}$ such that the map from $U_{u_{0}}$ to
 $Z_{s, T'}^{1/2}$, that associates to the initial data the unique
 solution of the IVP \eqref{ivp}, is  continuous, and in fact, real
 analytic.
    \end{theorem}

\begin{remark}
\label{scaling}
If $u(x,t)$ solves \eqref{ivp} then $u_\sigma (x,t) = \sigma^{1+a}
u(\sigma x , \sigma^{2+a} t)$ is also a solution, at least
formally. This dilation invariance distinguishes the scaling invariant
Sobolev index for \eqref{ivp} to be $s = \half -a$. We observe
therefore that Theorem \ref{main} is close to the scaling invariant
result near $a=0$ and assumes much more regularity on the data than
scaling suggests is required for $a$ near 1. Note also that the results
obtained here in the $a \geq \frac{4}{5}$ are inferior to those
obtained by Kenig, Ponce and Vega in \cite{KPV0} in the sense that we
require more spatial regularity of the data as well as integrability
against a weight.  Further remarks on the
near optimality of our results near $a = 0$ are made at the end of
Section 3.
\end{remark}

\begin{remark}
  \label{global?}
Theorem \ref{main} could be iterated to obtain a global well-posedness
result for \eqref{ivp} for initial data in the space $F^{s_*}$ if one
could prove that ${{\| u(t) \|}_{F^{s_*}}}$ is finite for all
time. Recalling \eqref{fs}, \eqref{h} and \eqref{enorm}, all that
remains to be shown is the finiteness for all time of ${{\| x u(t)
    \|}_{H^{-s_*}}}$. It may be possible that a Gronwall type argument
like that presented in Theorem 3.3 of \cite{SautKP} for the KP
equation will prove this. Note that the results in \cite{KPV0} already
provide global well-posedness for \eqref{ivp} in the energy space
identified using \eqref{h} in the regime $a \geq \frac{4}{5}.$

\end{remark}

The paper is organized as follows: in Section \ref{linear}
we present some estimates regarding the solution of the linear
homogeneous and  inhomogeneous IVP
associated with \eqref{ivp}. In Section \ref{bilinear} we present some
bilinear estimates needed in Section \ref{proof} where we prove
Theorem \ref{main} using a contraction method. The paper ends with
Section \ref{appendix}, an appendix in which we prove  some
estimates involving the spaces $X^{b}_{s}$ and
$Y^{b}_{s_{0},s_{1}}$ and cut-off functions in time.

\section{The linear estimates}\label{linear}
We denote with $W(t)u_{0}$ the solution of the linear IVP
\begin{equation}
\left\{ \begin{array}{l}
\partial_t u + \partial_xD_{x}^{1+a} u = 0, \\
 u(x,0) = u_0(x) \hspace{1.5cm} x \in \dbR, \,  t \in \dbR.
\end{array}\right.
\label{linivp}\end{equation}
The first lemma we present contains a priori estimates for the solution
$W(t)u_{0}$.

\begin{lemma}
For any $(\theta,\beta)\in [0,1]\times[0,a/2]$, we have:
\begin{eqnarray}
\label{hom}\left(\int_{\dbR}\|D^{\theta\beta/2}
    W(t)u_{0}\|^{q}_{p}\,dt\right)^{1/q}&\leq& C \|u_{0}\|_{L^{2}}\\
 \label{inhom}\left(\int_{\dbR}\|D^{\theta\beta}
    W(t-t')f(\cdot,t')\,dt'\|^{q}_{L^{p}_{x}}\,dt\right)^{1/q}&\leq& C
    \|f\|_{L^{q'}_{t}L^{p'}_{x}},
\end{eqnarray}
where $(q,p)=(2(2+a)/(\theta(\beta+1)), 2/(1-\theta))$ and
$1/p+1/p'=1/q+1/q'=1$.
    \end{lemma}
The proof of this lemma is due to Kenig, Ponce and Vega \cite{KPV1}.

A consequence of this lemma is the following Strichartz inequality:
\begin{corollary}
Let\footnote{Notice that for $a\in [0,1]$ we have
$b_{0}\in (0,1/2)$. } $b_{0}=(3+a)/(4(2+a))$, then for any
$f\in X^{b_{0}}_{0}$
\begin{equation}
    \label{l4}\|f\|_{L^{4}}\leq C\|f\|_{X^{b_{0}}_{0}}.
    \end{equation}
\end{corollary}
\begin{proof}
We first show that
\begin{equation}\label{l2}
    \|f\|_{L^{2}}\leq \|f\|_{X^{0}_{0}}.
    \end{equation}
In fact for a sequence of numbers $c_{j}$, $\|c_{j}\|_{l^{2}}\leq
\|c_{j}\|_{l^{1}}$ and it's enough to write
$$\|f\|_{L^{2}}=\left(\sum_{j\geq 0}\int\chi_{j}(\lambda-\om(\x))
|\hat f|^{2}(\lambda,\xi)\,d\x\,d\lambda\right)^{1/2},$$
and set $c_{j}=\left(\int\chi_{j}(\lambda-\om(\x))
|\hat f|^{2}(\lambda,\xi)\,d\x\,d\lambda\right)^{1/2}$. Next we take 
$\beta=0$ and
$p=q$ in \eqref{hom}, hence $\theta=(2+a)/(3+a)$ and $p=q=2(3+a)$.
A standard argument based on a layer decomposition along the
translates of the surface $S=\{(\x,\om(\x))\}$ (see for example
\cite{CKS1})   gives
\begin{equation}
    \label{23a}\|f\|_{L^{2(3+a)}}\leq C\|f\|_{X^{1/2}_{0}}.
    \end{equation}
We now interpolate \eqref{l2} and \eqref{23a}. We write
$1/4=\theta/2+(1-\theta)/2(3+a)$, $b_{0}=(1-\theta)/2$, and after a
simple calculation it follows that $\theta=(1+a)/(2(2+a))$
and $b_{0}=(3+a)/(4(2+a))$.
    \end{proof}
We  proceed now to the estimate of the group $W(t)$ in the spaces
$X^{b}_{s}$ and $Y^{b}_{s_{0},s_{1}}$. In the following we will always
assume that $\psi(t)$ is a smooth cut-off function supported in  the interval
$[-1,1]$.
\begin{lemma}\label{lem1}
For any $\delta\in (0,1)$ and $u_{0}\in H^{s}, \, s\in \dbR$, we have
$$\|\psi(t/\delta)W(t)u_{0}\|_{X^{1/2}_{s}}\leq C\|u_{0}\|_{H^{s}},$$
where $C>0$ is independent of $\delta$.
    \end{lemma}
\begin{proof}
Observe that
\begin{equation}\label{fwu0}
    \F(\psi(t/\delta)W(t)u_{0})(\x,\lambda)=\delta\hat{\psi}
(\delta(\lambda-\om(\x))\widehat{u_{0}}(\x).
\end{equation}
We thus need to estimate
$$\sum_{j\geq 0}2^{j/2}\left(\int_{\dbR^{2}}(1+|\x|)^{2s}
|\widehat{u_{0}}(\x)|^{2}\delta^{2}|\hat{\psi}(\delta(\lambda-\om(\x)))
|^{2}\chi_{j}(\lambda-\om(\x))\,d\x\,d\lambda\right)^{1/2}.$$
We first integrate the $\lambda$ variable, so we need to evaluate
$$I_{j}=\delta\left(\int_{\dbR}|\hat\psi(\delta\lambda)|^{2}
\chi_{j}(\lambda)\,d\lambda\right)^{1/2}.$$
For $j=0$
$$I_{0}=\delta\left(\int_{|\lambda|\leq 1}|\hat\psi(\delta\lambda)|^{2}
\,d\lambda\right)^{1/2}=\delta\left(\int_{|r|\leq \delta}
|\hat\psi(r)|^{2}\frac{dr}{\delta}\right)^{1/2}=C\delta.$$
For $j>0$
\begin{eqnarray*}
    I_{j}&=&\delta\left(\int_{2^{j-1}\leq |\lambda|\leq 2^{j}}
|\hat\psi(\delta\lambda)|^{2}\,d\lambda\right)^{1/2}
=\delta^{1/2}\left(\int_{\delta 2^{j-1}\leq |r|\leq \delta 2^{j}}
|\hat\psi(r)|^{2}\,dr\right)^{1/2}\\
&\leq&\delta^{1/2}\left(\int_{\delta 2^{j-1}\leq |r|\leq \delta 2^{j}}
|\hat\psi(r)|^{2}(1+|r|)^{2N}\frac{dr}{(1+|r|)^{2N}}\right)^{1/2}\\
&=&C\|\hat\psi(r)(1+|r|)^{N}\|_{L^{\infty}}\delta^{1/2}
\frac{(\delta 2^{j})^{1/2}}{(1+\delta 2^{j})^{N}},
\end{eqnarray*}
for any $N\in \dbN, N>1$. We then obtain
\begin{eqnarray*}
\|\psi(t/\delta)W(t)u_{0}\|_{X^{1/2}_{s}}&\leq&
\delta\|u_{0}\|_{H^{s}}\\
&+&\left(\sum_{j\geq 1}
\frac{\delta 2^{j}}{(1+\delta 2^{j})^{N}}\right)
\|\hat\psi(r)|(1+|r|)^{N}\|_{L^{\infty}}\|u_{0}\|_{H^{s}}.
\end{eqnarray*}
But it may be shown that $\sum_{j\geq 1}
\frac{\delta 2^{j}}{(1+\delta 2^{j})^{N}}\leq C$, uniformly for
$N>1$. Hence the lemma follows.
\end{proof}

The companion of the estimate we just proved is in the following lemma.
\begin{lemma}\label{lem7}
For any $\delta\in (0,1)$ and $u_{0}\in H^{s}, \, s\in \dbR$, we have
$$\|\psi(t/\delta)W(t)u_{0}\|_{Y^{1/2}_{s-2s_{*},s}}\leq C
\|u_{0}\|_{F^{s}},$$
where $C>0$ is independent of $\delta$.
    \end{lemma}
\begin{proof}
Using again \eqref{fwu0} we write
$$\frac{\partial}{\partial {\la}}
\F(\psi(t/\delta)W(t)u_{0})(\x,\lambda)=
\delta(\delta\hat{\psi}'
(\delta(\lambda-\om(\x))\widehat{u_{0}}(\x),$$
so that by Lemma \ref{lem1}
$$\|t\psi(t/\delta)W(t)u_{0}\|_{X^{1/2}_{s}}\leq
C\delta\|u_{0}\|_{H^{s}}.$$
On the other hand
\begin{eqnarray*}
    \frac{\partial}{\partial {\x}}
\F(\psi(t/\delta)W(t)u_{0})(\x,\lambda)&=&
\delta(\hat{\psi}
(\delta(\lambda-\om(\x))\frac{\partial}
{\partial {\x}}\widehat{u_{0}}(\x)\\
&-&\delta(\delta\hat{\psi}'
(\delta(\lambda-\om(\x))\om'(\x)\widehat{u_{0}}(\x).
\end{eqnarray*}
We recall that $|\om'(\x)|\sim |\x|^{1+a}$, hence by Lemma \ref{lem1}
$$\|x\psi(t/\delta)W(t)u_{0}\|_{X^{1/2}_{s-2s_{*}}}\leq C
\|xu_{0}\|_{H^{s-2s_{*}}}+C\delta\|u_{0}\|_{H^{s}},$$
as desired.
\end{proof}

We now need a version of Lemma \ref{lem1} and \ref{lem7} for the
solution of the inhomogeneous linear problem.

\begin{lemma}\label{lem8}
For any $\epsilon\in (0,1)$ and for any $\delta \in (0,1)$
$$\left\|\psi(t/\delta)\int_{0}^{t}W(t-t')h(t')\,dt'
\right\|_{X^{1/2}_{s}}\leq
C_{\epsilon}\delta^{-\epsilon}\|h\|_{X^{-1/2}_{s}}.$$
\end{lemma}
\begin{proof}
We follow the arguments of Kenig, Ponce and Vega in \cite{KPV3}. We write
\begin{eqnarray*}
& &\psi(t/\delta)\int_{0}^{t}W(t-t')h(t')\,dt'\\
&=& \psi(t/\delta)\int_{\dbR^{2}}e^{ix\x}\hat h(\x,\la)
\psi(\la-\om(\x))\frac{e^{i\la t}-e^{i\om(\x)
t}}{\la-\om(\x)}\,d\la\,d\x\\
&+&\psi(t/\delta)\int_{\dbR^{2}}e^{ix\x}\hat h(\x,\la)
[1-\psi(\la-\om(\x))]\frac{e^{i\la t}-e^{i\om(\x)
t}}{\la-\om(\x)}\,d\la\,d\x=I+II.
\end{eqnarray*}
To estimate $I$ we first perform a Taylor expansion
$$I=\sum_{k=1}^{\infty}\frac{i^{k}}{k!}t^{k}
\psi(t/\delta)\int_{\dbR}e^{ix\x+t\om(\x)}
\left(\int_{\dbR}\hat h(\x,\la)
\psi(\la-\om(\x))(\la-\om(\x))^{k-1}\,d\la\right)d\x.$$
Now let $t^{k}\psi(t/\delta)=\delta^{k}(t/\delta)^{k}\psi(t/\delta)=
\delta^{k}\psi_{k}(t/\delta)$, for $k\in \dbN$. Then
$$I=\sum_{k=1}^{\infty}\frac{i^{k}}{k!}\delta^{k}\psi_{k}(t/\delta)
W(t)G(x),$$
where
\begin{equation}\label{f}
\hat G(\x)=\int_{\dbR}\hat h(\x,\la)
\psi(\la-\om(\x))(\la-\om(\x))^{k-1}\,d\la.
\end{equation}
We want to use Lemma
\ref{lem1} together with its proof. To do so we first need
to estimate $|\widehat{\psi_{k}}(s)|$ and
$|\widehat{\psi_{k}}(s)|(1+|s|)^{N}$ for $N=2$, uniformly with
respect to $k\in \dbN$. Note that
$$|\widehat{\psi_{k}}(s)|\leq \int_{|t|\leq 1} |t|^{k}|\psi(t)|\,dt
\leq C,$$
uniformly in $k$. On the other hand
\begin{eqnarray*}
\widehat{\psi_{k}}(s)&=&\int_{|t|\leq 1}e^{ist}t^{k}\psi(t)\,dt\\
&=&(is)^{-2}\int_{|t|\leq 1}\partial_{t}^{2}(e^{ist})t^{k}\psi(t)\,dt
=(is)^{-2}\int_{|t|\leq 1}(e^{ist})\partial_{t}^{2}(t^{k}\psi(t))\,dt
\end{eqnarray*}
and for $|s|\geq 1$ it follows that
$$|\widehat{\psi_{k}}(s)|\leq C \frac{(1+k)^{2}}{(1+|s|)^{2}},$$
uniformly with respect to $k$. Then by Lemma \ref{lem1}
$$\|I\|_{X^{1/2}_{s}}\leq C \|G\|_{H^{s}},$$
and
\begin{eqnarray*}
\|G\|_{H^{s}}&=&\left(\int_{\dbR}(1+|\x|)^{2s}\left|\int_{\dbR}
\hat h(\x,\la)
\psi(\la-\om(\x))(\la-\om(\x))^{k-1}\,d\la\right|^{2}d\x\right)^{1/2}\\
&\leq&\left(\int_{\dbR}(1+|\x|)^{2s}\int_{|\la-\om(\x)|\leq 1}
|\hat h(\x,\la)|^{2}\,d\la d\x\right)^{1/2}\leq \|h\|_{X^{-1/2}_{s}},
    \end{eqnarray*}
as desired. This takes care of $I$. To estimate $II$ we write
$II=II_{1}+II_{2}$, where
\begin{eqnarray*}
II_{1}&=& -\psi(t/\delta)\int_{\dbR^{2}}e^{ix\x}e^{i\om(\x)t}
 [1-\psi(\la-\om(\x))]\frac{\hat h(\x,\la)}{\la-\om(\x)}\,d\la\,d\x,\\
II_{2}&=&\psi(t/\delta)\int_{\dbR^{2}}e^{ix\x}e^{i\la t}
    [1-\psi(\la-\om(\x))]\frac{\hat
    h(\x,\la)}{\la-\om(\x)}\,d\la\,d\x.
\end{eqnarray*}
In view of Lemma \ref{lem1}, to estimate $II_{1}$ we only need to show
that
\begin{equation}\label{ii1}
    \left\|\F^{-1}\left(\int_{\dbR}[1-\psi(\la-\om(\x))]
    \frac{\hat h(\x,\la)}{\la-\om(\x)}\,d\la\right)\right\|_{H^{s}}\leq C
    \|h\|_{X^{-1/2}_{s}}.
    \end{equation}
In fact
\begin{eqnarray*}
& 
&\left(\int_{\dbR}(1+|\x|)^{2s}\left(\int_{\dbR}[1-\psi(\la-\om(\x))]
\frac{\hat
    h(\x,\la)}{\la-\om(\x)}\,d\la\right)^{2}d\x\right)^{1/2}\\
&\leq&\left(\int_{\dbR}(1+|\x|)^{2s}\left(\sum_{j\geq 1}
\int_{\dbR}|\hat h(\x,\la)|\frac{\chi_{j}(\la-\om(\x))}{2^{j}}
\,d\la\right)^{2}d\x\right)^{1/2}\\
&\leq&\sum_{j\geq 1}\left(\int_{\dbR^{2}}(1+|\x|)^{2s}2^{-j}|\hat 
h(\x,\la)|^{2}
\chi_{j}(\la-\om(\x))\,d\la\,d\x\right)^{1/2},
    \end{eqnarray*}
    and \eqref{ii1} follows. Finally for $II_{2}$ we use Lemma
    \ref{lem5} and we reduce the matter to estimating
\begin{eqnarray*}
 & &\left\|\int_{\dbR^{2}}e^{ix\x+i\la t}
 [1-\psi(\la-\om(\x))]\frac{\hat h(\x,\la)}{\la-\om(\x)}\,d\la\,d\x
 \right\|_{X^{1/2}_{s}}\\
 &\leq&
 \sum_{j\geq 0}2^{j/2}\left(\int_{\dbR^{2}}
\frac{ |\hat h(\x,\la)|^{2}}{(1+|\la-\om(\x)|)^{2}}
(1+|\x|)^{2s}\chi_{j}(\la-\om(\x))\,d\la\,d\x\right)^{1/2}\\
&\leq&C\|h\|_{X^{-1/2}_{s}},
\end{eqnarray*}
and this concludes the proof of the lemma.
    \end{proof}
We now prove a corresponding lemma for the space $Y^{1/2}_{s-2s_{*},s}$.
\begin{lemma}\label{lem11}
For any $\epsilon\in (0,1)$ and for any $\delta \in (0,1)$
$$\left\|\psi(t/\delta)\int_{0}^{t}W(t-t')h(t')\,dt'
\right\|_{Y^{1/2}_{s-2s_{*},s}}\leq
C_{\epsilon}\delta^{-\epsilon}(\|h\|_{Y^{-1/2}_{s-2s_{*},s}}
+\|h\|_{X^{-1/2}_{s}}).$$
\end{lemma}
\begin{proof}
We perform a decomposition into $I$, $II_{1}$ and $II_{2}$
like in Lemma \ref{lem8}.
For $I$  we use Lemma \ref{lem7} and we only need to estimate
$$\|G\|_{F^{s}}=\|G\|_{H^{s}}+\|xG\|_{H^{s-2s_{*}}},$$
where $G$ is defined in \eqref{f}. In Lemma \ref{lem8} we proved
that $\|G\|_{H^{s}}\leq C\|h\|_{X^{-1/2}_{s}}$. To estimate
$\|xG\|_{H^{s-2s_{*}}}$ we first observe that
\begin{eqnarray*}
\frac{\partial}{\partial {\x}}\hat{G}(\x)&=&\int_{\dbR}
\frac{\partial}{\partial {\x}}\hat h(\x,\la)
\psi(\la-\om(\x))(\la-\om(\x))^{k-1}\,d\la\\
&-&(k-1)\om'(\x)\int_{\dbR}\hat h(\x,\la)
\psi(\la-\om(\x))(\la-\om(\x))^{k-2}\,d\la\\
&-&\om'(\x)\int_{\dbR}\hat h(\x,\la)
\psi'(\la-\om(\x))(\la-\om(\x))^{k-1}\,d\la=\widehat{J_{1}}+
\widehat{J_{2}}+\widehat{J_{3}}.
\end{eqnarray*}
We can now write
$$\|xG\|_{H^{s-2s_{*}}}\leq \sum_{i=1}^{3}\|J_{i}\|_{H^{s-2s_{*}}}.$$
Following the argument in the proof of Lemma \ref{lem8} it is easy to
see that
$$\|J_{1}\|_{H^{s-2s_{*}}}\leq C\|xh\|_{X^{-1/2}_{s-2s_{*}}}\leq
C\|h\|_{Y^{-1/2}_{s-2s_{*},s}}.$$
On the other hand, because $|\om'(\x)|\sim |\x|^{1+a}$, if follows that
$(1+|\x|)^{s-2s_{*}}|\om'(\x)| \lesssim (1+|\x|)^{s}$, hence
$$\sum_{i=1,2}\|J_{i}\|_{H^{s-2s_{*}}}\leq Ck \|h\|_{X^{-1/2}_{s}},$$
which suffices.

For $II_{1}$, we use again Lemma \ref{lem7} and we have to show that
\begin{equation}\label{ii1y}
    \left\|\F^{-1}\left(\int_{\dbR}[1-\psi(\la-\om(\x))]
    \frac{\hat h(\x,\la)}{\la-\om(\x)}\,d\la\right)\right
    \|_{F^{s}}\leq C
    (\|h\|_{Y^{-1/2}_{s-2s_{*},s}}+\|h\|_{X^{-1/2}_{s}}).
    \end{equation}
We observe that
\begin{eqnarray*}
& &\frac{\partial}{\partial {\x}}\left(\int_{\dbR}[1-\psi(\la-\om(\x))]
    \frac{\hat h(\x,\la)}{\la-\om(\x)}\,d\la\right)\\
&=& \int_{\dbR}\frac{\partial}{\partial {\x}}\hat h(\x,\la)
    \frac{[1-\psi(\la-\om(\x))]}{\la-\om(\x)}\,d\la \\
& +&\om'(\x)\int_{\dbR}\hat h(\x,\la)
    \frac{\psi'(\la-\om(\x))}{\la-\om(\x)}\,d\la
+\om'(\x)\int_{\dbR}\hat h(\x,\la)
\frac{[1-\psi(\la-\om(\x))]}{(\la-\om(\x))^{2}}\,d\la\\
&=&\widehat{M_{1}}+\widehat{M_{2}}+\widehat{M_{3}},
\end{eqnarray*}
and using both the argument in the proof of Lemma \ref{lem8} and the
observations above we have
\begin{eqnarray*}
    \|M_{1}\|_{H^{s}}&\leq& C\|xh\|_{X^{-1/2}_{s-2s_{*}}}\leq
C\|h\|_{Y^{-1/2}_{s-2s_{*},s}},\\
\sum_{i=2,3}\|M_{i}\|_{H^{s-2s_{*}}}&\leq& C \|h\|_{X^{-1/2}_{s}},
\end{eqnarray*}
which combined with \eqref{ii1} gives \eqref{ii1y}. Finally we turn to
$II_{2}$. We use Lemma \ref{lem9} and we obtain
\begin{eqnarray*}
 & &\left\|\psi(t/\delta)\int_{\dbR^{2}}e^{ix\x+i\la t}
 [1-\psi(\la-\om(\x))]\frac{\hat h(\x,\la)}{\la-\om(\x)}\,d\la\,d\x
 \right\|_{Y^{1/2}_{s-2s_{*},s}}\\
 &\leq& C_{\epsilon}\delta^{-\epsilon}\left
 \|\int_{\dbR^{2}}e^{ix\x+i\la t}
 [1-\psi(\la-\om(\x))]\frac{\hat h(\x,\la)}{\la-\om(\x)}\,d\la\,d\x
 \right\|_{Y^{1/2}_{s-2s_{*},s}}.
\end{eqnarray*}
We have that
\begin{eqnarray*}
& &\frac{\partial}{\partial {\x}}\left(\hat h(\x,\la)
    \frac{[1-\psi(\la-\om(\x))]}{\la-\om(\x)}\right)
=\frac{\partial}{\partial {\x}}\hat h(\x,\la)
    \frac{[1-\psi(\la-\om(\x))]}{\la-\om(\x)} \\
&+&\om'(\x)\hat h(\x,\la)
    \frac{\psi'(\la-\om(\x))}{\la-\om(\x)}
+\om'(\x)\hat h(\x,\la)
\frac{[1-\psi(\la-\om(\x))]}{(\la-\om(\x))^{2}},
\end{eqnarray*}
so we obtain
$$\left\|x\int_{\dbR^{2}}e^{ix\x+i\la t}
 [1-\psi(\la-\om(\x))]\frac{\hat h(\x,\la)}{\la-\om(\x)}\,d\la\,d\x
 \right\|_{X^{1/2}_{s-2s_{*}}}\leq C(\|h\|_{Y^{-1/2}_{s-2s_{*},s}}+
 \|h\|_{X^{-1/2}_{s}}).$$
 Next we turn to $\partial_{\la}$:
\begin{eqnarray*}
& &\frac{\partial}{\partial {\la}}\left(\hat h(\x,\la)
    \frac{[1-\psi(\la-\om(\x))]}{\la-\om(\x)}\right)=
    \frac{\partial}{\partial {\la}}\hat h(\x,\la)
    \frac{[1-\psi(\la-\om(\x))]}{\la-\om(\x)} \\
&-&\hat h(\x,\la)
    \frac{\psi'(\la-\om(\x))}{\la-\om(\x)}
+\hat h(\x,\la)
\frac{[1-\psi(\la-\om(\x))]}{(\la-\om(\x))^{2}}=\sum_{i=1,2,3}
\widehat{W_{i}},
\end{eqnarray*}
and following the proof in Lemma \ref{lem8} it is easy to see that
\begin{eqnarray*}
    \|W_{1}\|_{X^{1/2}_{s}}&\leq& C\|th\|_{X^{-1/2}_{s}}\\
    \sum_{i=2,3}\|(W_{i})\|_{X^{1/2}_{s}}&\leq &
    C\|h\|_{X^{-1/2}_{s}}.
    \end{eqnarray*}
The proof of the lemma is now complete.
    \end{proof}

\begin{remark} The results in the appendix, Lemma \ref{lem8} and Lemma
\ref{lem11} are probably also true with $\epsilon=0$.
    \end{remark}


\section{The bilinear estimates}\label{bilinear}
This section represent the heart of the matter of this paper. We start
with two propositions in which we prove bilinear estimates in the
spaces $X^{b}_{s_{*}}$ and $Y^{-b}_{-s_{*},s_{*}}$. After the proofs
of these estimates, we analyze in detail certain low/high frequency
interactions which reveal the near optimality of our bilinear estimates.
\begin{proposition}\label{prop1}
If $a\in (0,1)$ and $s_{*}=1/2+a/2$, there exists $0<b<1/2$ such that
\begin{eqnarray}\label{bil1}
    \|\partial_{x}(uv)\|_{X^{-b}_{s_{*}}}&\leq& C
    \|u\|_{X^{1/2}_{s_{*}}}\left(\|v\|_{X^{1/2}_{s_{*}}}
    +\|v\|_{X^{1/2}_{s_{*}}}^{1/2}
    \|v\|_{Y^{1/2}_{-s_{*},s_{*}}}^{1/2}\right)\\
\nonumber    &+&C\|v\|_{X^{1/2}_{s_{*}}}\left(\|u\|_{X^{1/2}_{s_{*}}}
    +\|u\|_{X^{1/2}_{s_{*}}}^{1/2}\|u\|_{Y^{1/2}_{-s_{*},s_{*}}}^{1/2}
    \right).
    \end{eqnarray}
    \end{proposition}
The companion estimate in the space $Y^{-b}_{-s_{*},s_{*}}$ takes the
following form:
\begin{proposition}\label{prop2}
If $a\in (0,1)$ and $s_{*}=1/2+a/2$, there exists $0<b<1/2$ such that
\begin{eqnarray}\label{bil2}
    \|\partial_{x}(uv)\|_{Y^{-b}_{-s_{*},s_{*}}}&\leq& C
    \|u\|_{Y^{1/2}_{-s_{*},s_{*}}}\left(\|v\|_{X^{1/2}_{s_{*}}}
    +\|v\|_{X^{1/2}_{s_{*}}}^{1/2}
    \|v\|_{Y^{1/2}_{-s_{*},s_{*}}}^{1/2}\right)\\
\nonumber    &+&C\|v\|_{Y^{1/2}_{-s_{*},s_{*}}}
\left(\|u\|_{X^{1/2}_{s_{*}}}
    +\|u\|_{X^{1/2}_{s_{*}}}^{1/2}\|u\|_{Y^{1/2}_{-s_{*},s_{*}}}^{1/2}
    \right)\\
 \nonumber    &+& C\|u\|_{X^{1/2}_{s_{*}}}\|v\|_{X^{1/2}_{s_{*}}}.
    \end{eqnarray}
    \end{proposition}
\noindent

\begin{proof}[Proof of Proposition \ref{prop1}]
We write
$$|\F(\partial_{x}(uv))(\x,\la)|=|\x|\left|
\int_{\dbR^{2}}\hat u(\xi_{1},\la_{1})
\hat v(\x-\xi_{1},\la-\la_{1})\,d\xi_{1}\,d\la_{1}\right|$$
and we split the set of integration  into the two regions
where $|\x-\xi_{1}|\leq |\x_{1}|$ and $|\x-\xi_{1}|>|\x_{1}|$. When we
integrate in the
first region we will obtain the first half of the right hand side of
\eqref{bil1}, when we integrate on the second region we will obtain the
second half. So without loss of generality we can assume that
$|\x-\xi_{1}|\leq |\x_{1}|$. Using duality and
setting $\x_{2}=\x-\xi_{1}$ allows us to rewrite the left hand
side of \eqref{bil1} as
\begin{eqnarray}\label{bil1lhs}
& &\sup_{\|g_{j}\|_{L^{2}}\leq 1}
\sum_{j\geq 0}\sum_{j_{1},j_{2}\geq 0}2^{-bj}\int_{|\x_{2}|\leq |\x_{1}|}
g_{j}(\x_{1}+\x_{2},\la_{1}+\la_{2})\chi_{j}(\la_{1}+\la_{2}
-\om(\x_{1}+\x_{2}))\\
\nonumber&\times&\frac{|\x_{1}+\x_{2}|\max(1,|\x_{1}+\x_{2}|)^{s_{*}}}
{\max(1,|\x_{1}|)^{s_{*}}\max(1,|\x_{2}|)^{s_{*}}}
\phi_{1,j_{1}}(\x_{1},\la_{1})\phi_{2,j_{2}}(\x_{2},\la_{2})
\,d\x_{1}\,d\la_{1}\,d\x_{2}\,d\la_{2},
\end{eqnarray}
where, if we set $u=u_{1}$ and $v=u_{2}$, then
$$\phi_{i,j_{i}}(\x_{i},\la_{i})=
|\hat u_{i}(\x_{i},\la_{i})|\chi_{j_{i}}(\la_{i}-\om(\x_{i}))
\max(1,|\x_{i}|)^{s_{*}},$$
and
$$\|u_{i}\|_{X^{1/2}_{s_{*}}}=\sum_{j_{i}\geq 0}2^{j_{1}/2}
\|\phi_{i,j_{i}}\|_{L^{2}}.$$
The analysis of \eqref{bil1lhs} is obtained by considering
different cases.

\noindent
{\bf Case 1:} $|\x_{1}|\leq 1$.\\
Then $|\x_{2}|\leq 1$ and $|\x_{2}+\x_{1}|\leq 2$. We observe that if
$G_{j}=\F^{-1}(g_{j}(\x,\la)\chi_{j}(\la-\om(\x)))$, then
$\|G_{j}\|_{L^{2}}\leq 1$. We use the
Strichartz inequality \eqref{l4} to obtain the bound
\begin{eqnarray*}
\eqref{bil1lhs} &\leq& \sum_{j\geq 0}\sum_{j_{i},j_{2}\geq 0}
2^{-jb}\|G_{j}\|_{L^{2}}\|\F^{-1}(\phi_{1,j_{1}})\|_{L^{4}}
\|\F^{-1}(\phi_{2,j_{2}})\|_{L^{4}}\\
&\leq&C\sum_{j\geq 0}\sum_{j_{i},j_{2}\geq 0}
2^{-jb}2^{j_{1}b_{0}}\|\phi_{1,j_{1}}\|_{L^{2}}
2^{j_{2}b_{0}}\|\phi_{2,j_{2}}\|_{L^{2}},
    \end{eqnarray*}
and this gives the desired estimate as long as $b>0$.

\noindent
{\bf Case 2:} $|\x_{1}+\x_{2}|\leq 1$ and $|\x_{1}|\geq 1$. \\
This case is treated like Case 1.

\noindent
{\bf Case 3:} $|\x_{1}|\geq  1, \, |\x_{1}+\x_{2}|\geq  1,$
and $1/4|\x_{1}|\leq |\x_{2}|\leq |\x_{1}|.$\\
Then by simple arguments
\begin{equation}\label{case31}
 \eqref{bil1lhs} \leq  \sum_{j\geq 0}\sum_{j_{i},j_{2}\geq 0}
2^{-jb}\int_{\dbR^{2}}g_{j}\chi_{j}|\x_{1}+\x_{2}|^{1-s_{*}}
\phi_{1,j_{1}}\phi_{2,j_{2}}\,d\x_{1}\,d\la_{1}\,d\x_{2}\,d\la_{2}.
\end{equation}
We next perform a dyadic decomposition of $\x_{1}$, hence $\x_{2}$, by
setting $\x_{i}\sim 2^{m_{1}}, \, m_{1}\sim m_{2}$. We then continue with
\begin{equation}\label{case32}
\eqref{case31}\leq  \sum_{m_{1}\geq 0}
\sum_{j\geq 0}\sum_{j_{i},j_{2}\geq 0}
2^{-jb}\int_{\dbR^{2}}g_{j}\chi_{j}|\x_{1}+\x_{2}|^{1-s_{*}}
\phi_{1,m_{1},j_{1}}\phi_{2,m_{1},j_{2}}\,d\x_{1}\,d\la_{1}
\,d\x_{2}\,d\la_{2}
    \end{equation}
where we used the notation
$$\phi_{i,m_{i},j_{i}}=\phi_{i,j_{i}}
\chi_{\{|\x_{i}|\sim 2^{m_{i}}\}}.$$
We have to consider two subcases

\noindent
{\bf Case 3a:} $j\geq (1-s_{*})m_{1}/b$. \\
We use again the Strichartz inequality \eqref{l4} to obtain the
following bound
\begin{equation}\label{case33}
\eqref{case32}\leq  \sum_{m_{1}\geq 0}\sum_{j\geq 
(1-s_{*})m_{1}/b}\sum_{j_{i},j_{2}\geq 0}
2^{-jb}2^{(1-s_{*})m_{1}}2^{j_{1}b_{0}}2^{j_{2}b_{0}}
\|\phi_{1,m_{1},j_{1}}\|_{L^{2}}\|\phi_{2,m_{1},j_{2}}\|_{L^{2}}.
    \end{equation}
We then sum in $j$ to get
$$
\eqref{case33}\leq  \sum_{m_{1}\geq 0}
\sum_{j_{i},j_{2}\geq 0}
2^{j_{1}b_{0}}2^{j_{2}b_{0}}
\|\phi_{1,m_{1},j_{1}}\|_{L^{2}}\|\phi_{2,m_{1},j_{2}}\|_{L^{2}},
$$
and  Cauchy-Schwarz in $m_{1}$ concludes the argument also in this
case.

\noindent
{\bf Case 3b:} $j\leq (1-s_{*})m_{1}/b$. \\
We change variables by setting
\begin{equation}\label{chan1}
    \la_{i}=\theta_{i}+\om(\x_{i}).\end{equation}
    Then we can continue the inequality in \eqref{case31} with
\begin{eqnarray}
    \label{case34}&\leq&  \sum_{m_{1}\geq 0}
\sum_{0\leq j\leq (1-s_{*})m_{1}/b}\sum_{j_{i},j_{2}\geq 0}
2^{-jb}\int_{\dbR^{2}}|\x_{1}+\x_{2}|^{1-s_{*}}g_{j}
(\x_{1}+\x_{2},\theta_{1}+\om(\x_{1})+\theta_{2}+\om(\x_{2}))\\
\nonumber&\times&
\chi_{j}(\theta_{1}+\om(\x_{1})+\theta_{2}+\om(\x_{2})-\om(\x_{1}+\x_{2}))
\phi_{1,m_{1},j_{1}}\phi_{2,m_{1},j_{2}}\,d\x_{1}\,d\x_{2}\,d\theta_{1}
\,d\theta_{2}.
    \end{eqnarray}
Note now that
\begin{equation}\label{relat}
    |\om(\x_{1})+\om(\x_{2})-\om(\x_{1}+\x_{2})|\leq
    C2^{\max(j_{1}, j_{2}, j)}.
    \end{equation}
At this point it is important to understand this expression in order
to obtain extra constraints that will allow a summation in the
expression given by \eqref{case34}.

We recall that we are analyzing the region
where $1/4|\x_{1}|\leq |\x_{2}|\leq |\x_{1}|$. Assume  in addition
that $\x_{1}\geq 0$ and $\x_{2}\geq 0$. Then we can write
$\x_{2}=\beta\xi_{1}$, where $1/4\leq \beta\leq 1$. Then
$$\om(\x_{1}+\x_{2})-\om(\x_{1})-\om(\x_{2})=
(\x_{1}+\x_{2})^{2+a}-\x_{1}^{2+a}-\x_{2}^{2+a}
=((1+\beta)^{2+a}-1-\beta^{2+a})\x_{1}^{2+a}.$$
We now consider the function
$$f(\beta)=(1+\beta)^{2+a}-1-\beta^{2+a},$$
and we observe that
$$f'(\beta)=(2+a)((1+\beta)^{1+a}-\beta^{1+a})\geq 0$$
and because $f(1/4)=c>0$ it follows that
$$|\om(\x_{1})+\om(\x_{2})-\om(\x_{1}+\x_{2})|\geq c|\x_{1}|^{2+a}.$$

If $\x_{1}\leq 0$ and $\x_{2}\leq 0$ the argument can be reapplied
thanks to the presence of the absolute value.

If $\x_{1}\geq  0$ and $\x_{2}\leq 0,$ then $\xi_{1}+\x_{2}\geq 0$.
If we set $\x_{2}=-\beta\x_{1}$, then again $1/4\leq \beta\leq 1$ and
$$\om(\x_{1}+\x_{2})-\om(\x_{1})-\om(\x_{2})=
(\x_{1}+\x_{2})^{2+a}-\x_{1}^{2+a}+|\x_{2}|^{2+a}
=((1-\beta)^{2+a}-1+\beta^{2+a})\x_{1}^{2+a}.$$
In this case we define the function
$$f(\beta)=1-(1-\beta)^{2+a}-\beta^{2+a},$$
and
$$f'(\beta)=(2+a)((1-\beta)^{1+a}-\beta^{1+a}).$$
We first analyze $f$ in  the range $1/2\leq \beta\leq 1$.
Here $f'(\beta)\leq 0$ and $f(1)=0$, so that $f(\beta)\geq 0$ as
$\beta\rightarrow 1$, and by Taylor expansion
$f(\beta)=-(2+a)(\beta-1)+O((\beta-1)^{2})$ and $f(\beta)\geq
c(1-\beta)$ for $1/2\leq \beta\leq 1$. Next, for
$1/4\leq \beta\leq 1/2$, $f'(\beta)\geq 0, \, f(1/4)>0$ hence also in
this case $f(\beta)\geq c(1-\beta)$, and so
$$|\om(\x_{1})+\om(\x_{2})-\om(\x_{1}+\x_{2})|\geq
c|\x_{1}|^{2+a}(1-\beta)=c|\x_{1}|^{1+a}|\x_{1}+\x_{2}|.$$
Finally we consider the case $\x_{1}< 0$ and $\x_{2}\geq 0.$ Then
$\x_{1}+\x_{2}\leq 0$ and we have
$$\om(\x_{1}+\x_{2})-\om(\x_{1})-\om(\x_{2})=
-|\x_{1}+\x_{2}|^{2+a}+|\x_{1}|^{2+a}-|\x_{2}|^{2+a},$$
which behaves like in the previous case. We summarize our findings as
follow: \\
if $1/4|\x_{1}|\leq |\x_{2}|\leq |\x_{1}|$ then
\begin{equation}\label{estrelat}
|\om(\x_{1})+\om(\x_{2})-\om(\x_{1}+\x_{2})|\geq
\left\{\begin{array}{l}
 c|\x_{1}|^{2+a} \, \, \mbox{ if } \, \, \x_{1}\x_{2}\geq  0\\
c|\x_{1}|^{1+a}|\x_{1}+\x_{2}| \, \, \mbox{ if } \, \,
\x_{1}\x_{2}\leq  0
\end{array}\right.
\end{equation}
Assume now that $j=\max(j,j_{1},j_{2})$. Then, using  \eqref{relat} and
\eqref{estrelat} we have that
$$2^{(1+a)m_{1}}\leq C2^{j}\leq C2^{(1-s_{*})m_{1}/b},$$
which is not possible since $(1+a)>(1-s_{*})/b$ if we take
$b$ such that $\frac{1-a}{2(1+a)}<b<1/2$. Thus, $j<\max(j_{1},j_{2})$.
Assume $j_{1}\geq j_{2}$, the other case is identical. We then
use our lower bound in \eqref{estrelat} to conclude that
$$|\x_{1}+\x_{2}|^{1-s_{*}}\leq
C2^{j_{1}(1-s_{*})}2^{-m_{1}(1+a)(1-s_{*})}.$$
We insert this estimate in  \eqref{case34} and we use the Strichartz
estimate \eqref{l4} to obtain
$$
\eqref{case34}\leq  \sum_{m_{1}\geq 0}
\sum_{j_{i},j_{2}\geq 0}\sum_{0\leq j\leq (1-s_{*})m_{1}/b}
2^{-jb}2^{jb_{0}}2^{j_{1}(1-s_{*})}2^{-m_{1}(1+a)(1-s_{*})}2^{j_{2}b_{0}}
\|\phi_{1,m_{1},j_{1}}\|_{L^{2}}\|\phi_{2,m_{1},j_{2}}\|_{L^{2}}.
$$
We now take $b>b_{0}, b<1/2$ and do the sum in $j$. Because
$(1-s_{*})=1/2-a/2<1/2$,  we  sum in $m_{1}$ and we obtain the desired
result.

\noindent
{\bf Case 4:} $|\x_{1}|\geq  1, \, |\x_{1}+\x_{2}|\geq  1,$
and $|\x_{2}|\leq 1/4|\x_{1}|.$\\
In this case $|\x_{1}+\x_{2}|\sim |\x_{1}| $ and
\begin{equation}\label{case41}
    \eqref{bil1lhs}\leq \sum_{j\geq 0}\sum_{j_{i},j_{2}\geq 0}
2^{-jb}\int_{\dbR^{2}}g_{j}\chi_{j}
\frac{|\x_{1}|}{\max(1,|\x_{2}|)^{s_{*}}}
\phi_{1,j_{1}}\phi_{2,j_{2}}\,d\x_{1}\,d\la_{1}\,d\x_{2}\,d\la_{2}.
\end{equation}
We first use the change of variables \eqref{chan1},
then we perform  a dyadic decomposition so that
$|\xi_{1}|\sim 2^{m_{1}}$
(and hence  $|\x_{1}+\x_{2}|\sim 2^{m_{1}}$), and we can continue
with
\begin{equation}\label{case42}
    \eqref{case41}\leq \sum_{m_{1}\geq 1}
    \sum_{j\geq 0}\sum_{j_{i},j_{2}\geq 0}
2^{-jb}\int_{\dbR^{2}}g_{j,m_{1}}\chi_{j}2^{m_{1}}
\phi_{1,m_{1},j_{1}}\phi_{2,j_{2}}
\,d\x_{1}\,d\x_{2}\,d\theta_{1}\,d\theta_{2}.
\end{equation}
We analyze two subcases.\\
{\bf Case 4a:} $m_{1}\leq (2b/(1-a) - \epsilon)j$, for some
$\epsilon>0$ to be chosen later.\\
Using an  argument similar to the one used in \cite{CKS1}, we
make the following change of variables:
\begin{equation}\label{chan2}
    (u,w)=T(\x_{1},\x_{2}),
    \end{equation}
where
\begin{eqnarray*}
    u &=&T_{1}(\x_{1},\x_{2})=\x_{1}+\x_{2}\\
   w &=&T_{2}(\x_{1},\x_{2})=
   \theta_{1}+\theta_{2}+\om(\x_{1})+\om(\x_{2}).
    \end{eqnarray*}
It is easy to see that if $J$ is the Jacobian of this change of
variables, then
\begin{equation}\label{jac}
    |J|=|\om'(\x_{1})-\om'(\x_{2})|\sim |\x_{1}|^{1+a}.
    \end{equation}
Now we define
$$H(\theta_{1},\theta_{2},u,w)=
\chi_{j}\phi_{1,m_{1},j_{1}}\phi_{2,j_{2}}
\circ T^{-1}(\theta_{1},\theta_{2},u,w),$$
and we write
\begin{eqnarray*}
    \eqref{case42}&\leq& \sum_{j\geq 0}
    \sum_{0\leq m_{1}\leq j(2b/(1-a) - \epsilon)}
    \sum_{j_{i},j_{2}\geq 0}
2^{-jb}2^{m_{1}}\int_{\dbR^{2}}g_{j,m_{1}}(u,w)
\frac{H(\theta_{1},\theta_{2},u,w)}{|J|}\,du\,dw\,d\theta_{1}
\,d\theta_{2}\\
&\leq&\sum_{j\geq 0}\sum_{0\leq m_{1}\leq j(2b/(1-a) - \epsilon)}
    \sum_{j_{i},j_{2}\geq 0}
    2^{-jb}2^{m_{1}}2^{-m_{1}(1+a)/2}\|g_{j,m_{1}}\|_{L^{2}}\\
  &\times&  \int_{|\theta_{i}|\sim 2^{j_{i}}}\left(\int\int
    \frac{H(\theta_{1},\theta_{2},u,w)^{2}}{|J|} \,du\,dw
    \right)^{1/2}\,d\theta_{1}\,d\theta_{2}.
\end{eqnarray*}
Now we invert the change of variable in the last integration in $u$
and $w$ to continue with
$$\leq \sum_{j\geq 0}\sum_{0\leq m_{1}\leq j(2b/(1-a) - \epsilon) }
    \sum_{j_{i},j_{2}\geq 0}
    2^{-jb}2^{m_{1}(1/2-a/2)}2^{j_{1}/2}2^{j_{2}/2}
\|g_{j,m_{1}}\|_{L^{2}}\|\phi_{1,m_{1},j_{1}}\|_{L^{2}}
\|\phi_{2,j_{2}}\|_{L^{2}},
$$
and if we sum with respect to $m_{1}$ we obtain
$$\leq
    \sum_{j\geq 0}\sum_{j_{i},j_{2}\geq 0}
    2^{-jb}2^{(2b/(1-a) - \epsilon)j(1/2-a/2)}
    2^{j_{1}/2}2^{j_{2}/2}\|\phi_{1,,j_{1}}\|_{L^{2}}
\|\phi_{2,j_{2}}\|_{L^{2}},
$$
and we finish by summing with respect to $j$, since
$-b+\left(2b/(1-a)-\epsilon\right)
\left(1/2-a/2\right)<0$. \\
{\bf Case 4b:} $m_{1}> j(2b/(1-a) - \epsilon)$, for some
$\epsilon>0$ to be chosen later.\\
This is the most delicate part of the proof. We return to
\eqref{bil1lhs} and we keep the function $v$ in it. Then we use the
change of variables \eqref{chan1} and we write
\begin{equation}\label{case4b1}
\eqref{bil1lhs}\leq \sum_{j\geq 0} \sum_{m_{1}\geq j(2b/(1-a) - \epsilon) }
    \sum_{j_{i},j_{2}\geq 0}2^{-jb}2^{m_{1}}
    \int g_{j,m_{1}}\chi_{j}\phi_{1,m_{1},j_{1}}
    |\hat v|\chi_{j_{2}}\,d\x_{1}\,d\x_{2}\,d\theta_{1}\,d\theta_{2}.
    \end{equation}
Next  for fixed $\x_{1},\, \theta_{1}$ and $\theta_{2}$,
we estimate the measure of the set $\Delta_{\x_{2}}$ such that
$|\theta_{1}+\theta_{2}+\om(\x_{1})+\om(\x_{2})-
\om(\x_{1}+\x_{2}))|\sim 2^{j}$. Thus, let
$$f(\x_{2})=\theta_{1}+\theta_{2}+\om(\x_{1})+\om(\x_{2})-
\om(\x_{1}+\x_{2}),$$
then
$$f'(\x_{2})=\om'(\x_{2})-\om'(\x_{1}+\x_{2})=(a+2)(|\x_{2}|^{1+a}
-|\x_{1}+\x_{2}|^{1+a}).$$
But in the region that we are analyzing
$ |\x_{2}|\leq 1/4|\x_{1}|$ and $|\x_{1}+\x_{2}|\geq
3/4|\x_{1}|$, hence
$$|f'(\x_{2})|\geq C|\x_{1}|^{1+a},$$
and by mean value theorem the measure of
$\Delta_{\x_{2}}$ can be estimated by
$$|\Delta_{\x_{2}}|\leq C2^{j}2^{-m_{1}(1+a)}.$$
We then use the change of variables \eqref{chan2} and we obtain
\begin{eqnarray}\label{case4b2}
\eqref{case4b1}&\leq& \sum_{j\geq 0}\sum_{m_{1}\geq j(2b/(1-a) - \epsilon) }
     \sum_{j_{i},j_{2}\geq 0}2^{-jb}2^{m_{1}(1-a)/2}
    2^{j_{1}/2}2^{j_{2}/2}2^{(j-m_{1}(1+a))/2}\\
\nonumber    &\times&
    \|g_{j,m_{1}}\|_{L^{2}}\|\phi_{1,m_{1},j_{1}}\|_{L^{2}}
    \||\hat v|\chi_{j_{2}}\|_{L^{2}_{\theta_{2}},L^{\infty}_{\x_{2}}}\\
\nonumber&\leq& \sum_{j\geq 0}\sum_{m_{1}\geq j(2b/(1-a) - \epsilon) }
     \sum_{j_{i},j_{2}\geq 0}2^{j(1/2-b)}
    2^{-m_{1}a}2^{j_{1}/2}2^{j_{2}/2}\\
\nonumber    &\times&\|g_{j,m_{1}}\|_{L^{2}}\|\phi_{1,m_{1},j_{1}}\|_{L^{2}}
    \||\hat v|(\x_{2},\theta_{2}+\om(\x_{2}))\chi_{j_{2}}
    (\theta_{2})\|_{L^{2}_{\theta_{2}},L^{\infty}_{\x_{2}}}.
    \end{eqnarray}
We now sum in $m_{1}$ to obtain
\begin{eqnarray}\label{case4b3}
\eqref{case4b2}&\leq& \sum_{j\geq 0} \sum_{j_{i},j_{2}\geq 0}
2^{j(1/2-b)}2^{-aj(2b/(1-a)-\epsilon)}2^{j_{1}/2}2^{j_{2}/2}\\
\nonumber&\times&
\|g_{j}\|_{L^{2}}\|\phi_{1,j_{1}}\|_{L^{2}}
    \||\hat v|(\x_{2},\theta_{2}+\om(\x_{2}))\chi_{j_{2}}
    (\theta_{2})\|_{L^{2}_{\theta_{2}},L^{\infty}_{\x_{2}}}.
\end{eqnarray}
Now observe that there exists $\epsilon>0$ such that
\begin{equation}\label{eb}
    (1/2-b)-a(2b/(1-a)-\epsilon)<0
    \end{equation}
if and only if
$(1/2-b)-a2b/(1-a)<0$ if and only if
$$b>\frac{1-a}{2(1+a)}$$
and  because for $0<a<1$ it follows that
$\frac{1-a}{2(1+a)}<1/2$, hence  there exists $b$ such that
$\frac{1-a}{2(1+a)}<b<1/2$ and $\epsilon$ such that \eqref{eb} is
satisfied. We go back to \eqref{case4b3} to sum in $j$ and we obtain
\begin{equation}\label{case4b4}
\eqref{case4b3}\leq  \sum_{j_{i},j_{2}\geq 0}
2^{j_{1}/2}2^{j_{2}/2}
\|\phi_{1,j_{1}}\|_{L^{2}}
    \||\hat v|(\x_{2},\theta_{2}+\om(\x_{2}))\chi_{j_{2}}
    (\theta_{2})\|_{L^{2}_{\theta_{2}},L^{\infty}_{\x_{2}}}.
\end{equation}
Next we note that if we set
$$|\hat v|(\x_{2},\theta_{2}+\om(\x_{2}))\chi_{j_{2}}(\theta_{2})=
|\hat v_{j_{2}}|(\x_{2},\theta_{2}+\om(\x_{2})),$$
then
\begin{eqnarray*}
|\hat v_{j_{2}}(\x_{2},\theta_{2}+\om(\x_{2}))|^{2}&\leq &
\int_{-\infty}^{\x_{2}}|\hat v_{j_{2}}(\nu_{2},\theta_{2}+
\om(\nu_{2}))|\left|\frac{\partial}{\partial {\nu_{2}}}
\hat v_{j_{2}}(\nu_{2},\theta_{2}+\om(\nu_{2}))\right|\,d\nu_{2}\\
&\leq &
\int_{-\infty}^{\x_{2}}|\hat v_{j_{2}}(\nu_{2},\theta_{2}+
\om(\nu_{2}))|(1+|\nu_{2}|)^{s_{*}}\frac{\left|\frac{\partial}{\partial {\nu_{2}}}
\hat v_{j_{2}}(\nu_{2},\theta_{2}+\om(\nu_{2}))\right|}
{(1+|\nu_{2}|)^{s_{*}}}\,d \nu_{2}\\
&\leq &\left(\int_{-\infty}^{\infty}|\hat v_{j_{2}}(\nu_{2},\theta_{2}+
\om(\nu_{2}))|^{2}(1+|\nu_{2}|)^{2s_{*}}\,d\nu_{2}\right)^{1/2}\\
&\times&\left(\int_{-\infty}^{\infty}
\left|\frac{\partial}{\partial {\nu_{2}}}
\hat v_{j_{2}}(\nu_{2},\theta_{2}+\om(\nu_{2}))\right|^{2}
\frac{d\nu_{2}}{(1+|\nu_{2}|)^{2s_{*}}}\right)^{1/2}.
\end{eqnarray*}
%
\begin{remark}
\label{LPonFTsideRemark} We pause to make a technical remark regarding
    the use of the spatial weight in the present analysis.
Observe that the $L^\infty_{\xi_2}$ norm first appeared in
\eqref{case4b2} and was just estimated using a Sobolev-type
inequality. The appearance of the differentiation operator
$\partial_{\nu_2}$ with respect to the Fourier variable corresponds
with the weight $x$ appearing in the definition \eqref{fs} of the space $F^s$
and the associated spacetime space $Y$ defined in \eqref{ys0s1b}.
\end{remark}

Because this estimate is independent of $\x_{2}$ we can write
\begin{eqnarray*}
\|\hat v_{j_{2}}(\x_{2},\theta_{2}+\om(\x_{2}))
\|_{L^{2}_{\theta_{2}}L^{\infty}_{\xi_{2}}}&\leq &
\left(\int_{\dbR}\left(\int_{-\infty}^{\infty}|\hat 
v_{j_{2}}(\nu_{2},\theta_{2}+
\om(\nu_{2}))|^{2}(1+|\nu_{2}|)^{2s_{*}}\,d\nu_{2}\right)^{1/2}
\right.\\
&\times &\left.\left(\int_{-\infty}^{\infty}
\left|\frac{\partial}{\partial {\nu_{2}}}
\hat v_{j_{2}}(\nu_{2},\theta_{2}+\om(\nu_{2}))\right|^{2}
\frac{d\nu_{2}}{(1+|\nu_{2}|)^{2s_{*}}}\right)^{1/2}\,d\theta_{2}
\right)^{1/2},
\end{eqnarray*}
and after Cauchy-Schwarz in $\theta_{2}$ we can continue with
\begin{eqnarray*}
&\leq&\left(\int_{\dbR^{2}}|\hat v_{j_{2}}(\nu_{2},\theta_{2}+
\om(\nu_{2}))|^{2}(1+|\nu_{2}|)^{2s_{*}}\,d\nu_{2}
\,d\theta_{2}\right)^{1/4}\\
&\times &\left(\int_{\dbR^{2}}
\left|\frac{\partial}{\partial {\nu_{2}}}
\hat v_{j_{2}}(\nu_{2},\theta_{2}+\om(\nu_{2}))\right|^{2}
\frac{1}{(1+|\nu_{2}|)^{2s_{*}}}\,d\nu_{2}\,d\theta_{2}
\right)^{1/4}.
\end{eqnarray*}
If we substitute this last estimate in \eqref{case4b4}, we sum in
$j_{1}$ and we use Cauchy-Schwarz in $j_{2}$, we obtain
\begin{eqnarray}\label{case4b5}
\eqref{case4b4}&\leq& \|u\|_{X^{1/2}_{s_{*}}}
\left(\sum_{j_{2}\geq 0} \|\hat v_{j_{2}}(\nu_{2},\theta_{2}+
\om(\nu_{2}))(1+|\nu_{2}|)^{s_{*}}\|_{L^{2}}2^{j_{2}/2}\right)^{1/2}\\
\nonumber&\times&\left(\sum_{j_{2}\geq 0} 
\left\|\frac{\partial}{\partial {\nu_{2}}}
\hat v_{j_{2}}(\nu_{2},\theta_{2}+\om(\nu_{2}))
(1+|\nu_{2}|)^{-s_{*}}\right\|_{L^{2}}2^{j_{2}/2}\right)^{1/2}.
    \end{eqnarray}
Now observe that
\begin{eqnarray*}
 \left|\frac{\partial}{\partial {\nu_{2}}}
\hat v_{j_{2}}(\nu_{2},\theta_{2}+\om(\nu_{2}))\right|&=&
\left|\frac{\partial\hat v_{j_{2}}}{\partial {\nu_{2}}}
(\nu_{2},\theta_{2}+\om(\nu_{2}))+
\om'(\nu_{2})\frac{\partial\hat v_{j_{2}}}{\partial {\la_{2}}}
(\nu_{2},\theta_{2}+\om(\nu_{2}))\right|\\
&\leq&\left|\frac{\partial\hat v_{j_{2}}}{\partial {\nu_{2}}}
(\nu_{2},\theta_{2}+\om(\nu_{2}))\right|+
C|\nu_{2}|^{1+a}\left|\frac{\partial\hat v_{j_{2}}}{\partial {\la_{2}}}
(\nu_{2},\theta_{2}+\om(\nu_{2}))\right|.
\end{eqnarray*}
Because $s_{*}=1/2+a/2$, it follows that
$$\frac{|\nu_{2}|^{1+a}}{(1+|\nu_{2}|)^{s_{*}}}\leq
|\nu_{2}|^{1/2+a/2}\leq (1+|\nu_{2}|)^{1/2+a/2},$$
and our desired inequality holds.
\end{proof}
We next turn to the proof of the corresponding fact in the $Y$ space.

\begin{proof}[Proof of Proposition \ref{prop2}]
We first estimate
\begin{equation}\label{dla}
    \frac{\partial}{\partial {\la}}\F(\partial_{x}
(uv))=\frac{\partial}{\partial {\la}}
\left(\xi\int_{\dbR}\hat u(\x_{1},\la_{1})
\hat v(\x-\x_{1},\la-\la_{1})\,d\x_{1}\,d\la_{1}\right).
\end{equation}
We again split the integration
into the region $|\x-\x_{1}|\leq |\x_{1}|$ and
$|\x-\x_{1}|\geq |\x_{1}|$. For symmetry reasons explained at the
beginning of the proof of Proposition \eqref{prop1}, it is enough to
estimate the integral on the first region. In this
case we replace \eqref{dla} with
$$
\frac{\partial}{\partial {\la}}\F(\partial_{x}
(uv))=\xi\int_{\dbR}\frac{\partial\hat u}{\partial {\la_{1}}}(\x_{1},\la_{1})
\hat v(\x-\x_{1},\la-\la_{1})\,d\x_{1}\,d\la_{1}.
$$
Since $\F^{-1}(\frac{\partial\hat u}{\partial {\la}})\in
X^{1/2}_{s_{*}}$, we obtain the desired estimate using the proof of
Proposition \ref{prop1}.

Next we turn to
\begin{equation}\label{dx}
    \frac{\partial}{\partial {\x}}\F(\partial_{x}
(uv))=\frac{\partial}{\partial {\x}}
(\xi (\hat u*\hat v)(\x,\la))=\hat u*\hat v(\x,\la)+\xi
\frac{\partial}{\partial {\x}}(\hat u*\hat v)(\x,\la).
\end{equation}
and we recall that we have to estimate it in $X^{-b}_{-s_{*}}$.
We start with  the term $\hat u*\hat v(\x,\la)$.
Our usual duality argument  leads us to rewriting the left hand
side of \eqref{bil2} as
\begin{eqnarray*}
\sup_{\|g_{j}\|_{L^{2}}\leq 1}
\sum_{j\geq 0}2^{-jb}\int_{\dbR^{4}}\frac{g_{j}(\x_{1}+\x_{2},
\la_{1}+\la_{2})\chi_{j}(\la_{1}+\la_{2}-\om(\x_{1}+\x_{2}))}
{\max(1,|\x_{1}+\x_{2}|)^{s_{*}}}\hat u(\x_{1},\la_{1})
\hat v(\x_{2},\la_{2})\,d\x_{1}\,d\x_{2}\,d\la_{1}\,d\la_{2},
    \end{eqnarray*}
and one can simply use the Strichartz inequality \eqref{l4}, just like
we did in Case 1 of  Proposition
\ref{prop1}. This leads to an estimate by
$C\|u\|_{X^{1/2}_{s_{*}}}^{1/2}\|v\|_{X^{1/2}_{s_{*}}}^{1/2}.$.

We now turn to the term $\xi
\frac{\partial}{\partial {\x}}(\hat u*\hat v)(\x,\la)$ In order to
treat this term we write
$$1=\theta_{1}(\rho)+\theta_{2}(\rho)+\theta_{3}(\rho),$$
where $supp \, \theta_{1}\subset \{|\rho|<1\}, \, \theta_{1}=1$ on
$\{|\rho|<1/2\}, \, supp \, \theta_{2}\subset \{1/2<|\rho|<2\}, \,
supp \, \theta_{3}\subset \{|\rho|>2\}, \, \theta_{3}=1$ on
$\{|\rho|>4\}$. Then
$$\xi \frac{\partial}{\partial {\x}}(\hat u*\hat v)(\x,\la)=
\sum_{i=1}^{3}\xi\frac{\partial}{\partial {\x}}
\int\hat u(\x_{1},\la_{1})\theta_{i}\left(\frac{|\x-\x_{1}|}{|\x_{1}|}
\right)\hat v(\x-\x_{1},\la-\la_{1})\,d\xi_{1}\,d\la_{1}
=I_{1}+I_{2}+I_{3}.$$
We start with $I_{1}$ by writing
\begin{eqnarray}
  \label{I1}  I_{1}&=&\xi
\int\hat u(\x_{1},\la_{1})\frac{\partial}{\partial {\x}}\left(
\theta_{1}\left(\frac{|\x-\x_{1}|}{|\x_{1}|}
\right)\right)\hat v(\x-\x_{1},\la-\la_{1})\,d\xi_{1}\,d\la_{1}\\
\nonumber&+&\x\int\hat u(\x_{1},\la_{1})
\theta_{1}\left(\frac{|\x-\x_{1}|}{|\x_{1}|}
\right)\frac{\partial}{\partial {\x}}
\hat v(\x-\x_{1},\la-\la_{1})\,d\xi_{1}\,d\la_{1}=I_{11}+I_{12},
\end{eqnarray}
Clearly
$$I_{11}=\xi
\int\hat u(\x_{1},\la_{1})\frac{sign \, (\x-\x_{1})}{|\x_{1}|}
\theta_{1}'\left(\frac{|\x-\x_{1}|}{|\x_{1}|}
\right)\hat v(\x-\x_{1},\la-\la_{1})\,d\xi_{1}\,d\la_{1}.$$
The contribution of $I_{11}$ leads us to estimating
\begin{eqnarray}
    \label{bil21}& &\sum_{j\geq 0}2^{-jb}
    \int_{\dbR^{4}}|g_{j}(\x_{1}+\x_{2},
\la_{1}+\la_{2})\chi_{j}(\la_{1}+\la_{2}-\om(\x_{1}+\x_{2}))|\\
\nonumber&\times&\frac{|\x_{1}+\x_{2}|}{\max(1,|\x_{1}+\x_{2}|)^{s_{*} 
}|\x_{1}|}
\left|\theta_{1}'\left(\frac{|\x_{2}|}{|\x_{1}|}
\right)\right||\hat u(\x_{1},\la_{1})|
|\hat v(\x_{2},\la_{2})|\,d\x_{1}\,d\x_{2}\,d\la_{1}\,d\la_{2},
    \end{eqnarray}
    and because in the support of $\theta_{1}'$ we have that
$|\x_{2}|\sim |\x_{1}|$ it follows that $|\x_{1}+\x_{2}|\leq
|\x_{1}|$, and hence \eqref{bil21} can be handled just by using the
Strichartz inequality \eqref{l4}.

In $I_{12}$ we can write
$$\frac{\partial}{\partial {\x}}
\hat v(\x-\x_{1},\la-\la_{1})=-\frac{\partial}{\partial {\x_{1}}}
\hat v(\x-\x_{1},\la-\la_{1})$$
and we  integrate by parts in $\x_{1}$ to get
\begin{eqnarray}
 \label{I12}   I_{12}&=&\x\int\frac{\partial}{\partial {\x_{1}}}
    \hat u(\x_{1},\la_{1})
\theta_{1}\left(\frac{|\x-\x_{1}|}{|\x_{1}|}
\right)
\hat v(\x-\x_{1},\la-\la_{1})\,d\xi_{1}\,d\la_{1}\\
\nonumber&+&\x\int\hat u(\x_{1},\la_{1})
\frac{\partial}{\partial {\x_{1}}}\theta_{1}
\left(\frac{|\x-\x_{1}|}{|\x_{1}|}
\right)
\hat v(\x-\x_{1},\la-\la_{1})\,d\xi_{1}\,d\la_{1}=I_{121}+I_{122}.
\end{eqnarray}
Now
$$\frac{\partial}{\partial {\x_{1}}}\theta_{1}
\left(\frac{|\x-\x_{1}|}{|\x_{1}|}
\right)=-\frac{sign \,(\x-\x_{1})}{|\x_{1}|}
\theta_{1}'\left(\frac{|\x-\x_{1}|}{|\x_{1}|}
\right)+\frac{sign \, \x_{1}|\x-\x_{1}|}{|\x_{1}|^{2}}
\theta_{1}'\left(\frac{|\x-\x_{1}|}{|\x_{1}|}
\right),$$
hence
$$\left|\frac{\partial}{\partial {\x_{1}}}\theta_{1}
\left(\frac{|\x-\x_{1}|}{|\x_{1}|}
\right)\right|\leq \frac{C}{|\x_{1}|}
\left|\theta_{1}'\left(\frac{|\x-\x_{1}|}{|\x_{1}|}
\right)\right|$$
and we go back to \eqref{bil21}. We turn our attention to $I_{121}$,
and notice that $|\x-\x_{1}|\leq |\x_{1}|$ in the region of
integration. We thus are led to estimate
\begin{eqnarray}
    \label{bil22}& &\sum_{j\geq 0}2^{-jb}
    \int_{|\x_{2}|\leq |\x_{1}|}g_{j}(\x_{1}+\x_{2},
\la_{1}+\la_{2})\chi_{j}(\la_{1}+\la_{2}-\om(\x_{1}+\x_{2}))
\\
\nonumber&\times&\frac{|\x_{1}+\x_{2}|}{\max(1,|\x_{1}+\x_{2}|)^{s_{*}}}
\frac{\max(1,|\x_{1}|)^{s_{*}}}{\max(1,|\x_{2}|)^{s_{*}}}
\frac{|\partial_{\x_{1}}\hat u(\x_{1},\la_{1})|}
{\max(1,|\x_{1}|)^{s_{*}}}
\max(1,|\x_{2}|)^{s_{*}}|\hat v(\x_{2},\la_{2})|
\,d\x_{1}\,d\x_{2}\,d\la_{1}\,d\la_{2}.
    \end{eqnarray}
We now consider several cases:

\noindent
{\bf Case 1:}$|\x_{1}|\leq 1$.\\
Then $|\x_{2}|\leq 1$ and $|\x_{1}+\x_{2}|\leq 2$ and we proceed just
like in Case 1 in the proof of Proposition \ref{prop1}.

\noindent
{\bf Case 2:}  $|\x_{1}|\geq 1, \, |\x_{1}+\x_{2}|\leq 1/2$.\\
Then $|\x_{2}|\sim |\x_{1}|$ and the multiplier in \eqref{bil22}
can be bounded by
$$\frac{|\x_{1}+\x_{2}|}{\max(1,|\x_{1}+\x_{2}|)^{s_{*}}}
\frac{\max(1,|\x_{1}|)^{s_{*}}}{\max(1,|\x_{2}|)^{s_{*}}}
\leq C,$$
and we are back in Case 1.

\noindent
{\bf Case 3:}  $|\x_{1}|\geq 1, \, |\x_{1}+\x_{2}|\geq 1/2$
 and $1/4|\x_{1}|\leq |\x_{2}|\leq |\x_{1}|$.\\
Then
$$\frac{|\x_{1}+\x_{2}|}{\max(1,|\x_{1}+\x_{2}|)^{s_{*}}}
\frac{\max(1,|\x_{1}|)^{s_{*}}}{\max(1,|\x_{2}|)^{s_{*}}}
\leq C|\x_{1}+\x_{2}|^{1-s_{*}},$$
and the estimate follows like in Case 3  of Proposition
\ref{prop1}.

\noindent
{\bf Case 4:}  $|\x_{1}|\geq 1, \, |\x_{1}+\x_{2}|\geq 1/2$
 and $1/4|\x_{1}|\geq |\x_{2}|$.\\
Then $|\x_{1}|\sim |\x_{1}+\x_{2}|$ and
$$\frac{|\x_{1}+\x_{2}|}{\max(1,|\x_{1}+\x_{2}|)^{s_{*}}}
\frac{\max(1,|\x_{1}|)^{s_{*}}}{\max(1,|\x_{2}|)^{s_{*}}}
\leq C\frac{|\x_{1}|}{\max(1,|\x_{2}|)^{s_{*}}},$$
and the result follows like in Case 4  of Proposition
\ref{prop1}. This concludes the estimate involving $I_{1}$.

We next turn to $I_{2}$. Just like  in \eqref{I1},  we can write
$$I_{2}=I_{21}+I_{22},$$
where now  both integrals  involve $\theta_{2}$. Because
$supp \, \theta_{2}'\subset\{1/2<\rho<2\}$, for $I_{21}$
we can repeat the argument presented for $I_{11}$. Now using the same
splitting described in \eqref{I12}, we can write
$I_{22}=I_{221}+I_{222}$, where again $\theta_{1}$ is replaced by
$\theta_{2}$. It is easy to see that $I_{222}$ can be treated just like
$I_{11}$. We are left with $I_{221}$, which leads to:
\begin{eqnarray}
    \label{bil23}& &\sum_{j\geq 0}2^{-jb}
    \int_{|\x_{2}|\leq |\x_{1}|}g_{j}(\x_{1}+\x_{2},
\la_{1}+\la_{2})\chi_{j}(\la_{1}+\la_{2}-\om(\x_{1}+\x_{2}))
\theta_{2}\left(\frac{|\x_{2}|}{|\x_{1}|}\right)\\
\nonumber&\times&
\frac{|\x_{1}+\x_{2}|}{\max(1,|\x_{1}+\x_{2}|)^{s_{*}}}
\frac{\max(1,|\x_{1}|)^{s_{*}}}{\max(1,|\x_{2}|)^{s_{*}}}
\frac{|\partial_{\x_{1}}\hat u(\x_{1},\la_{1})|}
{\max(1,|\x_{1}|)^{s_{*}}}
\max(1,|\x_{2}|)^{s_{*}}|\hat v(\x_{2},\la_{2})|
\,d\x_{1}\,d\x_{2}\,d\la_{1}\,d\la_{2}.
    \end{eqnarray}
We note that $|\x_{2}|\sim|\x_{1}|$ on $supp \, \theta_{2}$. So if
$|\x_{1}|\leq 1$, then $|\x_{2}|, |\x_{2}+\x_{1}| \leq C$ and the
bound follows like in Case 1 of Proposition \ref{prop1}.
If $|\x_{1}|\geq 1$ and $|\x_{2}+\x_{1}| \leq 1$, then again the bound
follows because
$$
\theta_{2}\left(\frac{|\x_{2}|}{|\x_{1}|}\right)
\frac{|\x_{1}+\x_{2}|}{\max(1,|\x_{1}+\x_{2}|)^{s_{*}}}
\frac{\max(1,|\x_{1}|)^{s_{*}}}{\max(1,|\x_{2}|)^{s_{*}}}\leq C,
$$
and we go back to the previous case. Finally,
if $|\x_{1}|\sim |\x_{2}|\geq 1$ and $|\x_{2}+\x_{1}| \geq 1$
we have
$$
\theta_{2}\left(\frac{|\x_{2}|}{|\x_{1}|}\right)
\frac{|\x_{1}+\x_{2}|}{\max(1,|\x_{1}+\x_{2}|)^{s_{*}}}
\frac{\max(1,|\x_{1}|)^{s_{*}}}{\max(1,|\x_{2}|)^{s_{*}}}
\leq C|\x_{2}+\x_{1}|^{1-s_{*}},
$$
and we proceed like in Case 3 of Proposition \ref{prop1}. Thus also
for $I_{2}$ we obtain the desired bound.

We finally turn to $I_{3}$, which we split in the usual way  into
$I_{3}=I_{31}+I_{33}$. Note that $supp \, \theta_{3}'
\subset\{2<\rho<4\}$, and so $I_{31}$ is handled like $I_{11}$. Now
observe that
$$I_{32}=\x\int\hat u(\x_{1},\la_{1})
\theta_{3}\left(\frac{|\x-\x_{1}|}{|\x_{1}|}
\right)\frac{\partial}{\partial_{\x}}
\hat v(\x-\x_{1},\la-\la_{1})\,d\xi_{1}\,d\la_{1}$$
and if we make the change of variables $\nu_{1}=\x-\x_{1}, \,
\tau_{1}=\la-\la_{1}$, we obtain
$$I_{32}=\x\int\hat u(\x-\nu_{1},\la-\tau_{1})
\theta_{3}\left(\frac{|\nu_{1}|}{|\x-\nu_{1}|}
\right)\frac{\partial}{\partial_{\nu_{1}}}
\hat v(\nu_{1},\tau_{1})\,d\nu_{1}\,d\tau_{1},$$
and on $supp \, \theta_{3}$ we have $|\nu_{1}|\geq 2|\x-\nu_{1}|
>|\x-\nu_{1}|$. We then treat this term as $I_{121}$, with the roles
of $v$ and $u$ interchanged.

\bigskip
\noindent{\bf{Analysis of certain low/high frequency interactions}}

We first show that if the $Y$-norms in \eqref{bil1} are ignored, the
resulting estimate failes. This analysis shows that a contraction
mapping argument in the space $X_s^b$ alone will fail, which is 
a special case of the negative result in \cite{MST}. Next, we
show, by analyzing a more refined low/high frequency interaction, that
the natural extension of \eqref{bil1} when $a=0$ fails to hold. This
demonstrates the near optimality of \eqref{bil1} near $a=0$.

The following proposition is contained in \cite{MST}.
\begin{proposition}
\label{XfailProp}
Consider the space $X_s^b$ adapted to the dispersive function
$\omega_a ( \xi) = |\xi |^{1+a} \xi $ defined in \eqref{xsb} for $a
\in [0,1)$. The bilinear estimate
\begin{equation}
  \label{NoY}
  {{\| \partial_x (u_1 u_2 ) \|}_{X_{s_*}^{-b}}} \leq C {{\| u_i
      \|}_{X_{s_*}^b}} {{\| u_j \|}_{X_{s_*}^b}} 
\end{equation}
where $(i,j) = (1,2)$ FAILS to hold for any $b \in \R$.
\end{proposition}
\begin{proof}
As in \cite{MST}, consider
\begin{equation*}
  {\widehat{u_1}} ( \xi , \lam ) = \chi_{[\alpha/2 , \alpha]} ( \xi )
  \chi_{\{|\lam - \omega_a ( \xi ) | \leq 1 \} } ( \lam ),
\end{equation*}
\begin{equation*}
  {\widehat{u_2}} ( \xi , \lam ) = \chi_{[N  , N +\alpha]} ( \xi )
  \chi_{\{|\lam - \omega_a ( \xi ) | \leq 1 \} } ( \lam ),
\end{equation*}
where $\chi_S$ denotes a smooted version of the characteristic
function of the set S and $\alpha \thicksim N^{-1 -a}$.

\begin{remark}
  \label{geometric}
The choice $\alpha \thicksim N^{-1 -a}$ may be explained
geometrically. Consider the Taylor expansion of the dispersive
function $\omega_a ( \xi)$ at $\xi = N$ for some $N \gg 1$:
\begin{equation*}
  \omega_a ( N + x ) = N^{2+a} + (2+a) N^{1+a} x + (2+a) (1+a) N^a x^2
  + \dots.
\end{equation*}
Note that the horizontal line $\lam = N^{2+a} $ is the zeroth order
approximation to the curve $\lam = \omega_a ( \xi )$ at the point $(N,
N^{2+a} )$. The curve $\lam = \omega (N + x) $ separates a distance
$\thicksim 1$ from the horizontal line 
$\lam = N^{2+a}$ when $x \thicksim N^{-1-a}$. A similar
geometric observation regarding the tangent line to the curve at $(N,
N^{2+a})$ explains the choice $\beta \thicksim N^{-a/2}$ below.
\end{remark}

By analyzing a convolution, it may be shown that 
\begin{equation}
  \label{left}
  {\widehat{[\partial_x (u_1 u_2)]}} ( \xi, \lam ) \thicksim N \alpha 
\chi_{\{ |\lam - \omega_a ( \xi) | \lesssim 1\} } (\lam ) \chi_{[N, N+
  \alpha ]}(\xi ).
\end{equation}
Thus, the left-side of \eqref{NoY} is of size $N^{s_*} N \alpha
\alpha^\half$. The right-side of \eqref{NoY} is of the size $N^{s_*}
\alpha$. For \eqref{NoY} to hold, we must have
\begin{equation*}
  N^{s_*} \alpha \alpha^\half \lesssim N^{s_*} \alpha
\end{equation*}
for all $N \gg 1$. But $\alpha \thicksim N^{-1 -a}$ so this requires
$\half - \frac{a}{2} \leq 0$ forcing $a \geq 1$. Therefore \eqref{NoY}
fails.
\end{proof}

Next, we consider a more refined low/high frequency interaction
between
\begin{equation*}
 {\widehat{u_1}} ( \xi , \lam ) = \chi_{[-\beta , \beta]} ( \xi )
  \chi_{\{|\lam - \omega_a ( \xi ) | \leq 1 \} } ( \lam ),
\end{equation*}
\begin{equation*}
  {\widehat{u_2}} ( \xi , \lam ) = \chi_{[N  , N +\beta]} ( \xi )
  \chi_{\{|\lam - \omega_a ( \xi ) | \leq 1 \} } ( \lam ),
\end{equation*} 
where $\beta = N^{-a/2}$. Recall that this choice of $\beta$ may be
explained geometrically by noting that when the tangent to the graph
of the dispersive function at $(N, N^{2+a})$ separates $\thicksim 1$
from the graph. Note that the support of $\widehat{u_1}$ intersects
lines $\lambda = constant$ in intervals of length at most $\beta =
N^{-\alpha/2}.$ Note also that the support of $\widehat{u_2}$
intersects lines $\lambda = constant$ in intervals of size at most
$N^{-1 -a}$ and that the projection of the support of $\widehat{u_2}$
along the $\lam$-axis is an interval of length $N^{1 + a/2}$. A
convolution analysis shows that 
\begin{equation}
  \label{lefts}
  {\widehat{[\partial_x (u_1 u_2)]}} ( \xi, \lam ) \thicksim N N^{-1-a} 
\chi_{\{ |\lam - N^{2+a} | \lesssim N^{1+a/2} \}} (\lam ) \chi_{\{|\xi - N|
  \lesssim \beta\}}(\xi ).
\end{equation}
We may now calculate the left and right sides of \eqref{bil1} and find
\begin{equation}
  \label{leftside}
  {\mbox{left-side of}}~\eqref{bil1} \gtrsim N^{s_*} N N^{-1 -a }
  \beta^\half \sum_{0 \leq j \leq 10^{-10} \log N } 2^{j ( -b)} 2^{j /2}.
\end{equation}
The right-side of \eqref{bil1} involves various terms of three sizes
\begin{equation*}
  N^{s_*} \beta, ~\beta N^{\frac{1+a}{2}} , ~ \beta N^{s_*} N^{\frac{a}{4}}.
\end{equation*}
The last of these is the largest since $s_* = \half + \frac{a}{2}$ and
$a \in [0,1)$. For \eqref{bil1} to hold, we must have
\begin{equation}
  \label{required}
  N^{s_*} N N^{-1 -a } \beta^\half \sum_{0 \leq j \leq 10^{-10} \log N
  } 2^{j ( -b)} 2^{j /2}
\lesssim \beta N^{s_*} N^{\frac{a}{4}},
\end{equation}
for all $N \gg 1.$ The sum contributes $\thicksim N^\epsilon $ if $ b = \half - \epsilon$
and contributes $\thicksim \log N$ if $b = \half$. Recall that $b \in
(\frac{1-a}{2(1+a)}, \half)$ was used in the proof of Proposition \ref{prop1}.
Since $\beta \thicksim N^{-a/2}$, \eqref{required} requires
\begin{equation*}
  N^{-a} \log N \lesssim 1
\end{equation*}
which fails to hold if $a = 0$. Note also that for $a>0$ there is a
bit of slack allowing us to take $b<\half$ for this interaction which
is consistent with Proposition \ref{prop1}.

We summarize the negative result just obtained with the following
statement.
\begin{proposition}
  \label{azfailure}
Consider the space $X_s^b$ adapted to the dispersive function
$\omega_a ( \xi) = |\xi |^{1+a} \xi $ defined in \eqref{xsb} for
$a=0$. 
The bilinear estimate
\begin{equation}
  \label{azfails}
  {{\| \partial_x (u_1 u_2 ) \|}_{X_{\half}^{-\half}}} \leq C {{\| u_i
      \|}_{X_{\half}^\half}}
\left( {{\| u_j \|}_{X_{\half}^\half}} +\| u_j \|_{X^{1/2}_{\half}}^{1/2}
    \|u_j \|_{Y^{1/2}_{-\half,\half}}^{1/2}\right)  + (i ~{\mbox{switched with}}~ j)
\end{equation}
where $(i,j) = (1,2)$ FAILS to hold.

\end{proposition}

\end{proof}


\section{Proof of Theorem \ref{main}}\label{proof}
We now turn to the local well-posedness result for the IVP \eqref{ivp}
stated in Theorem \ref{main}. With the results proved in Section
\ref{linear}, Section \ref{bilinear} and Section \ref{appendix},
the proof becomes standard,  similar to
the one presented in \cite{KPV3}. For completeness we
present it in  detail.

We first assume that $s=s_{*}$ and  that  $u_{0}\in F^{s_{*}}$ and we choose
$\psi\in C^{\infty}_{0}(\dbR)$, a cut-off function for the interval
$[-1,1]$. For $\delta$ small, to be chosen, we consider the nonlinear
mapping
$$\Phi_{u_{0}}(v)=\psi(t/\delta)W(t)u_{0}+
\psi(t/\delta)\int_{0}^{t}W(t-t')\partial_{x}(v^{2})dt'.$$
We then consider the ball
$$B_{a}=\{v\in Z^{1/2}_{s_{*}}/ \|v\|_{Z^{1/2}_{s_{*}}}\leq a\}.$$
Theorem \ref{main} follows from proving that there exist $a>0$ and
$\delta=\delta(\|u_{0}\|_{F^{s_{*}}})$ such that
\begin{equation}\label{ball}
   \Phi_{u_{0}}:B_{a}\longrightarrow B_{a},
    \end{equation}
and
\begin{equation}\label{cont}
   \|\Phi_{u_{0}}(v)- \Phi_{u_{0}}(w)\|_{Z^{1/2}_{s_{*}}}\leq
   1/2\|v-w\|_{Z^{1/2}_{s_{*}}}.
    \end{equation}
By the linear estimates in Lemma \ref{lem1} and Lemma \ref{lem7}
we have that
$$\|\psi(t/\delta)W(t)u_{0}\|_{Z^{1/2}_{s_{*}}}\leq
C\|u_{0}\|_{F^{s_{*}}}.$$
Choose $a=2C\|u_{0}\|_{F^{s_{8}}}$. We need to estimate
the nonlinear part of the mapping $\Phi_{u_{0}}$. We choose a new
cut-off function in time $\tilde\psi$ such that $\tilde\psi=1$ on the
support of $\psi$. Then
$$NL(t,x)=\psi(t/\delta)\int_{0}^{t}W(t-t')\partial_{x}(v^{2})dt'=
\psi(t/\delta)\int_{0}^{t}W(t-t')\partial_{x}(\tilde{v}^{2})dt',$$
where $\tilde{v}=\tilde\psi(t/\delta) v$. Thus by Lemma
\ref{lem8} and Lemma  \ref{lem11},
$$\|NL\|_{Z^{1/2}_{s_{*}}}\leq C_{\epsilon}\delta^{-\epsilon}
\|\partial_{x}(\tilde{v}^{2})\|_{Z^{-1/2}_{s_{*}}}
$$
and by Lemma \ref{lem12} we can continue with
$$\leq C_{\epsilon}\delta^{-\epsilon+\theta}
\|\partial_{x}(v^{2})\|_{Z^{-b}_{s_{*}}}$$
and by the bilinear estimates \eqref{bil1} and \eqref{bil2}
$$\leq C_{\epsilon}\delta^{-\epsilon+\theta}
\|v\|_{Z^{1/2}_{s_{*}}}^{2}\leq
C_{\epsilon}\delta^{-\epsilon+\theta}a^{2}.$$
We fix now $b$ as in the bilinear estimates \eqref{bil1}
and \eqref{bil2}, $\theta$ corresponding to $b$, and we
choose $\epsilon=\theta/2$. Then we choose $\delta$ small enough such
that $C_{\theta}\delta^{\theta/2}a<1/2$. We the obtain \eqref{ball}.
To prove \eqref{cont} we argue similarly, using
$\partial_{x}(v^{2})-\partial_{x}(w^{2})=\partial_{x}((v+w)(v-w))$
and the bilinear estimates to obtain
$$\|\Phi_{u_{0}}(v)- \Phi_{u_{0}}(w)\|_{Z^{1/2}_{s_{*}}}\leq
   2aC_{\theta}\delta^{\theta}\|u-w\|_{Z^{1/2}_{s_{*}}},$$
   and we now choose $\delta$ so small that
   $2C_{\theta}\delta^{\theta}a<1/2.$ This finishes the proof of the
theorem when $s=s_{*}$. If $s>s_{*}$ then we set $\tilde s=s-s_{*}$
and we use the following higher order bilinear estimate
\begin{eqnarray*}\|\partial_{x}v^{2}\|_{X^{-b}_{s}}&\leq&
\|\partial_{x}(vD^{\tilde s}v)\|_{X^{-b}_{s_{*}}}
\leq C\|v\|_{Z^{1/2}_{s_{*}}}\|v\|_{Z^{1/2}_{s}}
\end{eqnarray*}
and similarly for the space $Y^{-b}_{s-s_{*},s}$.
Now it is easy to see that by repeating the chain of inequalities
 for $s>s_{*}$, one can prove that also in this case the time interval $\delta$
depends only on $\|u_{0}\|_{F^{s_{*}}}$.


\section{Appendix: Time cutoffs and the $X$ and $Y$ norms }\label{appendix}
In this section we first present a few lemmas that quantify how the norm of
a function $f$ changes in the spaces $X^{1/2}_{s}$ and
$Y^{1/2}_{s-2s_{*},s}$ when multiplied by a time cut-off function
relative to an interval of size $\delta$.

A good estimate would be
\begin{equation}\label{better}
    \|\psi(t/\delta)f\|_{X^{1/2}_{s}}\leq C\|f\|_{X^{1/2}_{s}},
    \end{equation}
where $C$ is independent of $\delta$ and $s$. Such an estimate
is a bit cumbersome to prove. So, we replace \eqref{better}
with a weaker estimate, that in any case will suffice for what we need.
\begin{lemma}\label{lem5}
Given $\epsilon>0$, there exists $C_{\epsilon}>0$ such that
$$\|\psi(t/\delta)f\|_{X^{1/2}_{s}}\leq C_{\epsilon}
\delta^{-\epsilon}\|f\|_{X^{1/2}_{s}},$$
for any $s \in \dbR$.
\end{lemma}
If we are willing to relax a bit the regularity of the left hand
side in \eqref{better} then we have:
\begin{lemma}\label{lem6}
Given $b\in (0,1/2)$, there exists $\theta=\theta(b)>0$ such that
$$\|\psi(t/\delta)f\|_{X^{b}_{s}}\leq C\delta^{\theta}
\|f\|_{X^{1/2}_{s}},$$
for any $s \in \dbR$.
\end{lemma}
To prove  Lemma \ref{lem5} we need to introduce the  auxiliary space
$\tilde X^{b}_{s}$ defined as the closure of the Schwartz functions
with respect to the norm\footnote{One can also write
$$\|f\|_{\tilde X^{b}_{s}}=\left(\sum_{j\geq 0}2^{2jb}\int_{\dbR^{2}}
\chi_{j}(\la-\om(\x))(1+|\x|)^{2s}|\hat f(\x,\la)|^{2}d\x\,d\la
\right)^{1/2}$$.}
$$\|f\|_{\tilde X^{b}_{s}}=\left(\int_{\dbR^{2}}
(1+|\la-\om(\x)|)^{2b}(1+|\x|)^{2s}|\hat f(\x,\la)|^{2}d\x\,d\la
\right)^{1/2}.$$
Lemma \ref{lem5} will be proved as an interpolation between two
estimates in the space $\tilde X^{b}_{s}$, one with $b<1/2$ and one
with $b>1/2$. More precisely, we  use the following results:
\begin{lemma}\label{lem3}
If  $b\in (0,1/2)$, then
\begin{equation}\label{lem31}
    \|\psi(t/\delta)f\|_{\tilde X^{b}_{s}}\leq C
\|f\|_{\tilde X^{b}_{s}},\end{equation}
for any $s\in \dbR$.

If $b\in (1/2,1)$, then
\begin{equation}\label{lem32}
    \|\psi(t/\delta)f\|_{\tilde X^{b}_{s}}\leq C
\delta^{1/2-b}\|f\|_{\tilde X^{b}_{s}},\end{equation}
for any $s \in \dbR$.

In both cases $C$ is independent of $\delta$.
\end{lemma}
\begin{proof}
The proof of this lemma follows  the proof of Lemma 3.1 in
\cite{KPV3}. We first observe that using the notation
$$\widehat{J^{s}h}(\x)=(1+|\x|)^{s}\hat h(\x) \, \, \,
\mbox{ and } \, \, \, \, \widehat{\Lambda^{b}m}(\la)
=(1+|\la|)^{b}\hat m(\la)$$
we have
$$\|f\|_{\tilde X^{b}_{s}}\sim\|\Lambda^{b}W(t)J^{s}f\|_{L^{2}}.$$
We also note that for any $s$ and any $b\geq 0$
\begin{eqnarray*}
\left(\int_{\dbR^{2}}(1+|\x|)^{2s}|\widehat{\psi_{\delta} f}|^{2}
\,d\x\,d\lambda\right)^{1/2}
&=&\|\psi(\delta^{-1}\cdot) J^{s}f\|_{L^{2}}
\leq C \|J^{s}f\|_{L^{2}}\leq \|f\|_{\tilde X^{b}_{s}},
\end{eqnarray*}
where $\psi_{\delta}(t)=\psi(\delta^{-1}t)$.   If we also use the
notation
$$\widehat{D^{b}m}(\la)=|\la|^{b}\hat m(\la),$$
then
\begin{eqnarray*}
\int_{\dbR^{2}}|\lambda-\om(\x)|^{2b}(1+|\x|)^{2s}
|\widehat{ \psi_{\delta}f}|^{2}\,d\x\,d\lambda
&=&\int_{\dbR}(1+|\x|)^{2s}\int_{\dbR}|\lambda-\om(\x)|^{2b}
|\widehat{ \psi_{\delta} f}|^{2}\,d\lambda\,d\x\\
&=&\int_{\dbR}(1+|\x|)^{2s}\int_{\dbR}|\lambda|^{2b}
|\hat f*(\delta\hat\psi(\delta\cdot))(\lambda-\om(\xi))
|^{2}\,d\lambda\,d\x\\
&=&\int_{\dbR}(1+|\x|)^{2s}\left(\int_{-\infty}^{\infty}
|D^{b}(e^{i\om(\x)t}\tilde{f}(t)\psi(t/\delta))|^{2}\,dt\right)\,d\x,
\end{eqnarray*}
where in the last step we denoted with $\tilde{f}$ the partial Fourier
transform with respect to $\xi$, and we used the identity
\begin{equation}\label{identity}
\|D^{b}(e^{iat}f(t))\|_{L^{2}}^{2}=\int_{-\infty}^{\infty}
|\hat f(\lambda)|^{2}|\lambda-a|^{2b}d\lambda.
    \end{equation}
    We thus need to prove that for any $a\in \dbR$
\begin{equation}\label{lem33}
    \int_{-\infty}^{\infty}
|D^{b}(e^{iat}f(t)\psi(t/\delta))|^{2}\,dt\leq
C\int_{-\infty}^{\infty}|\hat f(\lambda)|^{2}|\lambda-a|^{2b}d\lambda.
    \end{equation}
We recall the Leibniz rule for fractional derivative \cite{KPV2}:
for any $\alpha\in(0,1)$
\begin{equation}\label{leib}
\|D_{t}^{\alpha}(fg)-gD^{\alpha}_{t}(f)-fD^{\alpha}_{t}(g)\|_{L^{2}}
\leq C \|D^{\alpha_{1}}_{t}(f)\|_{L^{p_{1}}}
\|D^{\alpha_{2}}_{t}(g)\|_{L^{p_{2}}},
    \end{equation}
where $\alpha=\alpha_{1}+\alpha_{2}, \, \alpha_{i}\geq 0$
and $1/2=1/p_{1}+1/p_{2}$.
In our case we use \eqref{leib} with $\alpha_{1}=b$ and $\alpha_{2}=0$.
Then
$$\|D^{b}(e^{iat}f(t)\psi(t/\delta))\|_{L^{2}}\leq C
\|D^{b}(e^{iat}f(t))\psi(t/\delta)\|_{L^{2}}+
\|e^{iat}f(t)D^{b}\psi(t/\delta)\|_{L^{2}}
=I_{1}+I_{2}.$$
To estimate $I_{1}$ we use the fact that $\|\psi\|_{L^{\infty}}\leq C$
and the identity \eqref{identity}.
We are then left with the term $I_{2}$. By the Sobolev embedding
theorem
$$\|e^{iat}f\|_{L^{q}}\leq C\|D^{b}(e^{iat}f)\|_{L^{2}},$$
provided $1/q=1/2-b$. We now set $ 2r=q$, it follows that
$1/r=1-2b>0$ and we can apply H\"older inequality to obtain
\begin{equation}\label{i2}
    I_{2}=\|e^{iat}f(t)D^{b}\psi(t/\delta)\|_{L^{2}}\leq
\|e^{iat}f\|_{L^{2r}}\|D^{b}\psi(t/\delta)\|_{L^{2r'}}.
\end{equation}
Now observe that $1/r'=1-1/r=2b$, hence
$$\|D^{b}\psi(t/\delta)\|_{L^{2r'}}=\delta^{1/2r'-b}
\|D^{b}\psi(t)\|_{L^{2r'}}=C.$$
It follows that
$$I_{2}\leq C\|D^{b}(e^{iat}f)\|_{L^{2}},$$
and \eqref{lem33} follows thanks to \eqref{identity}.

The proof of  \eqref{lem32}
is exactly the same as the one given to prove Lemma 3.1 in
\cite{KPV3}, hence we decided to omit it.
\end{proof}
\begin{proof}[Proof of Lemma \ref{lem5}]
We first observe the
following facts about real interpolation. Let $A$ be a Banach space,
and
$$l^{q}_{s}(A)=\left\{(f_{j}) / f_{j}\in A, \, \left(\sum_{j\geq 0}
(2^{js}\|f_{j}\|_{A})^{q}\right)^{1/q}<\infty\right\}.$$
Then by Theorem 5.6.1 in \cite{BL}
\begin{equation}\label{interpolation}
    l^{q}_{s_{\theta}}(A)=(l^{q_{0}}_{s_{0}}(A),
l^{q_{1}}_{s_{1}}(A))_{\theta,q}, \, \, \, (s_{0}\ne s_{1})
\end{equation}
 where
$s_{\theta}=\theta s_{0}+(1-\theta) s_{1}, \, 1\leq q\leq \infty,
\, 1\leq q_{i}\leq \infty, \, i=1,2$.
Using \eqref{interpolation} we see that for any $s\in \dbR$,
$$(\tilde X^{b_{1}}_{s}, \tilde X^{b_{2}}_{s},)_{\theta,1}=
X^{b}_{s},$$
where $b=\theta b_{1}+(1-\theta)b_{2}$, and $b_{1}\ne b_{2}$.
Now Lemma \ref{lem5} follows
from this observation,
interpolation and Lemma \ref{lem3}.
\end{proof}
\begin{proof}[Proof of Lemma \ref{lem6}]
Let us fix $b<b_{1}<b_{2}<1/2$. Assume one can prove
\begin{equation}\label{lem61}
    \|\psi(t/\delta)f\|_{\tilde X^{b_{1}}_{s}}\leq C\delta^{\theta}
\|f\|_{\tilde X^{b_{2}}_{s}}.
\end{equation}
Then, since $\|\psi(t/\delta)f\|_{ X^{b}_{s}}
\leq \|\psi(t/\delta)f\|_{\tilde X^{b_{1}}_{s}}$ and
$\|f\|_{\tilde X^{b_{2}}_{s}}\leq
\|f\|_{ X^{1/2}_{s}}$, we will be done. By complex interpolation,
\eqref{lem61} will follow from the two estimates below:

\begin{equation}\label{lem62}
    \|\psi(t/\delta)f\|_{\tilde X^{b_{2}}_{s}}\leq C
\|f\|_{\tilde X^{b_{2}}_{s}},
\end{equation}
and
\begin{equation}\label{lem63}
    \|\psi(t/\delta)f\|_{\tilde X^{0}_{s}}\leq C\delta^{b_{2}}
\|f\|_{\tilde X^{b_{2}}_{s}}.
\end{equation}
Clearly \eqref{lem62} is true by \eqref{lem31} of Lemma \ref{lem3}.
So we only have to show \eqref{lem63}. By repeating some of the
arguments in  the proof of Lemma \ref{lem3}, one can prove that
\begin{equation}\label{lem64}
    \|\psi(t/\delta)f\|_{\tilde X^{0}_{s}}^{2}\sim
\int_{\dbR}\left(\int_{-\infty}^{\infty}
|(e^{i\om(\x)t}J^{s}f(t)\psi(t/\delta))|^{2}\,dt\right)\,d\x.
\end{equation}
Now we observe that
$$\|e^{iat}J^{s}f(t)\psi(t/\delta)\|_{L^{2}}\leq
\|e^{iat}J^{s}f\|_{L^{2r}}\|\psi(t/\delta)\|_{L^{2r'}},$$
where $1/q=1/2r=1/2-b_{2}$. Then we can continue with
$$\leq C\|e^{iat}J^{s}f\|_{L^{2r}}\delta^{1/2r'}\leq
C\|D^{b_{2}}(e^{iat}J^{s}f)\|_{L^{2}}\delta^{b_{2}}.$$
If we insert this in \eqref{lem64} and
use the identity \eqref{identity},  we then obtain
\eqref{lem63}.

\end{proof}
Lemma \ref{lem5} and \ref{lem6} can
also be proved if one replaces the spaces $X^{b}_{s}$
with  the spaces $Y^{b}_{s-2s_{*},s}$.
\begin{lemma}\label{lem9}
Given $\epsilon>0$, there exists $C_{\epsilon}>0$ such that
$$\|\psi(t/\delta)f\|_{Y^{1/2}_{s-s_{*},s}}\leq C_{\epsilon}
\delta^{-\epsilon}\|f\|_{Y^{1/2}_{s-s_{*},s}},$$
for any $s \in \dbR$.
\end{lemma}
\begin{proof}
Let $\psi_{\delta}(t)=\psi(t/\delta)$, so that
$$\F(\psi(t/\delta)f)(\x,\la)=\widehat{\psi_{\delta}}*_{\la}\hat f
(\x,\la),$$
hence
\begin{eqnarray*}
    \frac{\partial}{\partial {\x}}
\F(\psi(t/\delta)f)(\x,\la)&=&\widehat{\psi_{\delta}}*_{\la}
\frac{\partial}{\partial {\x}}\hat f(\x,\la),\\
\frac{\partial}{\partial {\la}}
\F(\psi(t/\delta)f)(\x,\la)&=&\widehat{\psi_{\delta}}*_{\la}
\frac{\partial}{\partial {\la}}\hat f(\x,\la).
\end{eqnarray*}
Then Lemma \ref{lem5} gives
\begin{eqnarray*}
    \|x\psi(t/\delta)f\|_{X^{1/2}_{s-2s_{*}}}&=&
C_{\epsilon}\delta^{-\epsilon}\|xf\|_{X^{1/2}_{s-2s_{*}}}\\
\|t\psi(t/\delta)f\|_{X^{1/2}_{s}}&=&
C_{\epsilon}\delta^{-\epsilon}\|tf\|_{X^{1/2}_{s}},
\end{eqnarray*}
and this concludes the proof.
    \end{proof}
With a similar argument we can also prove, from Lemma \ref{lem6},
\begin{lemma}\label{lem10}
Given $b\in (0,1/2)$, there exists $\theta=\theta(b)>0$ such that
$$\|\psi(t/\delta)f\|_{Y^{b}_{s-s_{*},s}}\leq C\delta^{\theta}
\|f\|_{Y^{b}_{s-s_{*},s}},$$
for any $s \in \dbR$.
\end{lemma}
We conclude this section with the following lemma.
\begin{lemma}\label{lem12}
For any $b \in (0,1/2)$, there exists $\theta=\theta(b)>0$ such that
\footnote{The space $Z^{b}_{s}$ is defined in \eqref{zs}.}
$$\|\psi(t/\delta)v\|_{Z^{-1/2}_{s}}\leq \delta^{\theta}
\|v\|_{Z^{-b}_{s}},$$
for any $s\in \dbR$.
    \end{lemma}
\begin{proof}
We only estimate $\|\psi(t/\delta)v\|_{X^{-1/2}_{s}}$ because
the estimate for $\|\psi(t/\delta)v\|_{Y^{-1/2}_{s-s_{*},s}}$
will then follow like Lemma \ref{lem9} followed from Lemma \ref{lem5}.
It is easy to
see that if $C^{1/2}_{-s}$ is the space defined through the norm
$$\|g\|_{C^{1/2}_{-s}}=\sup_{j} 2^{j/2}\left(
\int_{\dbR^{2}}(1+|\x|)^{2s}|\hat g(\x,\la)\|^{2}
\chi_{j}(\la-\om(\x))\,d\la\,d\x\right)^{1/2},$$
then $C^{1/2}_{-s}=(X^{-1/2}_{s})^{*}$, with duality pairing
$\int_{\dbR^{2}}v\bar g\,d\la\,d\x$. Thus
$$\|\psi(t/\delta)v\|_{X^{-1/2}_{s}}=\sup_{\|g\|_{C^{1/2}_{-s}}\leq 1}
\left|\int_{\dbR^{2}}\psi(t/\delta)v\bar g\,d\la\,d\x\right|=
\sup_{\|g\|_{C^{1/2}_{-s}}\leq 1}
\left|\int_{\dbR^{2}}v\overline{\psi(t/\delta)g}\,d\la\,d\x\right|.$$
We thus need to prove that
$$\|\psi(t/\delta)g\|_{C^{b}_{-s}}\leq C\delta^{\theta}
\|g\|_{C^{1/2}_{-s}}.$$
We go back to the proof of Lemma \ref{lem6}, where we showed that
for $0<b_{1}<b_{2}<1/2$, we have
$$\|\psi(t/\delta)g\|_{\tilde X^{b_{1}}_{-s}}\leq C\delta^{\theta}
\|g\|_{\tilde X^{b_{2}}_{-s}}.$$
Then
\begin{eqnarray*}
    \|\psi(t/\delta)g\|_{C^{b}_{-s}}&\leq&
\|\psi(t/\delta)g\|_{\tilde X^{b_{1}}_{-s}}\\
&\leq& C\delta^{b_{2}}\|g\|_{\tilde X^{b_{2}}_{-s}}
\leq C\delta^{b_{2}}\|g\|_{C^{1/2}_{-s}},
\end{eqnarray*}
as desired.
    \end{proof}

\end{document}